\documentclass[11pt]{article}
\usepackage{amsmath,amsthm,amsfonts,amssymb, eucal, amscd}
\usepackage{mathrsfs}
\usepackage{color}

\newtheorem{thm}[equation]{Theorem}
\newtheorem{cor}[equation]{Corollary}
\newtheorem{lem}[equation]{Lemma}
\newtheorem{prop}[equation]{Proposition}
\theoremstyle{definition}
\newtheorem{defn}[equation]{Definition}

\newtheorem{rem}[equation]{Remark}
\newtheorem{exam}[equation]{Example}
\newtheorem{exams}[equation]{Examples}

\newtheorem{subsec}[equation]{} 
\numberwithin{equation}{section}

\newcommand{\ann}{\hbox{\rm Ann}}
\newcommand{\psito}{\stackrel {\Psi}{\longrightarrow}}

\newcommand{\kk}{\mathbb{F}}
\newcommand{\z}{\mathbb{Z}}

\newcommand{\setof}[2]{\{ #1 \, | \, #2 \} }

\begin{document}
\title{Whittaker Modules for  \\ 
Generalized Weyl Algebras}
\author{Georgia Benkart \  
and \ 
Matthew Ondrus
}

\date{March 25, 2008} 

\maketitle

\begin{center}\textbf {ABSTRACT}
\begin{itemize}
\item [] {\small We investigate Whittaker modules for generalized Weyl algebras, a class of associative algebras which includes the quantum plane, Weyl algebras, the universal enveloping algebra of $\mathfrak {sl}_2$ and of Heisenberg Lie algebras, Smith's generalizations of $U(\mathfrak{sl}_2)$, various quantum analogues of these algebras,  and  many others.  We show that the Whittaker modules $V = Aw$  of the generalized Weyl algebra $A = R(\phi,t)$ are in bijection with the $\phi$-stable left ideals of $R$.   We determine the annihilator $\ann_A(w)$ of the cyclic generator $w$ of $V$.  We also describe the annihilator ideal $\ann_A(V)$ under certain assumptions that hold for most of the examples mentioned above.  As one special case,  we recover Kostant's well-known results on Whittaker modules and their associated annihilators for $U(\mathfrak{sl}_2)$.}
\end{itemize}
\end{center}

\section{Introduction}\label{sec:notationGWA}
 
In this work we study the notion of a Whittaker module in the setting of generalized Weyl algebras.  Generalized Weyl algebras were introduced by  Bavula \cite{bavula:ottawa92} and have been studied  extensively since then (see for example, \cite{bavula:gwar93}, \cite{rosenberg:quantizedAlg}, \cite{drozd:wmgwa96}).  We shall use the definition of a generalized Weyl algebra (GWA) given in \cite[1.1]{bavula:gwar93}.

Suppose that $R$ is a unital associative algebra over a field  $\kk$ with $\phi = (\phi_i)_{i \in \mathcal I}$ a collection of pairwise commuting automorphisms of $R$ indexed by the set $\mathcal I$ (which may be finite or infinite), and let $t = (t_i)_{i \in \mathcal I}$ be a collection of nonzero central elements of $R$ also indexed by $\mathcal I$.    The \emph{generalized Weyl algebra} $A = R(\phi, t)$ with base ring $R$ is the associative algebra generated over $R$ by elements $X_i$ and $Y_i$  for $i \in \mathcal I$  with defining relations 
\begin{eqnarray}
\label{eqn:YiXi} Y_iX_i=t_i &\hspace{.4in}& X_iY_i= \phi_i(t_i) \\
\label{eqn:Xir} X_ir = \phi_i(r) X_i &\hspace{.4in}& Y_i\phi_i(r) = r Y_i
\end{eqnarray}
for $r \in R$,   and 
\begin{equation}\label{eqn:XiXj}
[X_i,X_j] = [Y_i,Y_j] = [X_i,Y_j] = 0
\end{equation}
for $i \neq j$, where $[\ ,\ ]$ denotes the commutator $[a,b] = ab-ba$. 
 
\textbf{We always assume that the algebra $R$ is a domain which is left Noetherian}.  Thus by \cite[Prop.~1.3]{bavula:gwar93}, the algebra $A = R(\phi, t)$ is a domain,  
and $A$ is left Noetherian if $\mathcal I$ is finite.   
This assumption forces the automorphisms $\phi_i$ to satisfy 
$\phi_i(t_j) = t_j$ for $j \neq i$,  which can be seen from the calculation,

$$t_j X_i = Y_jX_j X_i = X_i Y_j X_j = X_i t_j = \phi_i(t_j) X_i.$$

Weyl algebras provide the prototypical examples of generalized Weyl algebras.  Let  $R = \kk[t_i \mid i \in \mathcal I]$, the polynomial algebra over $\kk$ in commuting variables $t_i$,  and let $\phi_i$ for $i \in \mathcal I$  be the automorphism of $R$ defined by  $\phi_i (t_j) = t_j - \delta_{i,j}$.   Assume $A =  R(\phi,t)$, where  the relations in \eqref{eqn:YiXi}, \eqref{eqn:Xir}, and  \eqref{eqn:XiXj} hold for these choices.    Then  $[Y_i,X_j] =\delta_{i,j}$ for all $i,j \in \mathcal I$,  and  $A$ is   a Weyl algebra realized as a generalized Weyl algebra.  

To construct  a second family of examples,   let $R = \kk[c, t_i \mid i \in \mathcal I]$, the polynomial algebra over $\kk$ in commuting variables $c,t_i, i \in \mathcal I$.    Let $\phi_i$ be the automorphism given by  $\phi_i(t_j)= t_j -\delta_{i,j}c$ and $\phi_i(c) = c$.    In the generalized Weyl algebra  $A = R(\phi,t)$ constructed from this data, $[Y_i,X_j] =\delta_{i,j}c$, and $c$ is central in $A$.  Thus,  $A= R(\phi,t)$ is 
isomorphic to the universal enveloping algebra of a Heisenberg Lie algebra in this case.   
The Weyl and Heisenberg algebras are always generalized Weyl algebras,  but $\mathcal I$ needs to be finite for $R$ to be  Noetherian.
     
The notion of a generalized Weyl algebra encompasses many more examples such as the universal enveloping algebra $U(\mathfrak {sl}_2)$  and quantized enveloping algebras $U_q(\mathfrak{sl}_2)$, $U_{r,s}(\mathfrak {sl}_2)$  of the Lie algebra $\mathfrak {sl}_2$, the Noetherian down-up algebras of  \cite{ben:dup},  \cite{benroby:downup}, generalized Heisenberg algebras, and quantum Weyl algebras.  We will explain many of these examples later as we discuss results on  their Whittaker modules.

Kostant \cite{kostant:wvrt78} introduced a class of modules for finite-dimensional complex 
semisimple Lie algebras and called them Whittaker modules because of their connections with Whittaker equations in number theory.  These modules have been studied subsequently in a variety of
different settings.  Mili$\check{\rm c}$i$\acute{\rm c}$ and Soergel \cite{milSoergel:compSeries} investigated  modules for semisimple Lie algebras  induced from Whittaker modules for parabolic subalgebras.   Whittaker modules for semisimple Lie algebras also appeared  in
the work of Brundan and Kleshchev \cite{brunKlesh:Walg} on shifted Yangians and $W$-algebras. 
  Christodoulopoulou \cite{christo:thesis} used Whittaker modules
for Heisenberg Lie algebras to construct irreducible modules for affine Lie algebras. 

 In \cite{block:am81}, Block showed that the simple modules for $\mathfrak{sl}_2$ over $\mathbb C$ are either  highest (or lowest) weight modules, Whittaker modules, or modules obtained by localization.   Whittaker modules for $U_q(\mathfrak{sl}_2)$ were investigated in \cite{ondrus:thesis}, \cite{ondrus:Whittaker},  where many analogues of Kostant's results on annihilators for Whittaker modules were shown to hold.   Because of the prominent role that Whittaker modules play in the representation theory of $\mathfrak{sl}_2$ and of its quantum analogues, we were motivated to study them in the context of generalized Weyl algebras as a way of providing a unified approach to these modules.

Fix $R, \phi,$ $t,$ and $A = R(\phi, t)$ as above, and let $\zeta = (\zeta_i)_{i \in \mathcal I}$ be a set of nonzero elements of $\kk$ indexed by $\mathcal I$.    We say that an $A$-module $V$ is a {\it Whittaker module of type $\zeta$}   if there exists $w \in V$ such that 
\begin{enumerate} 
\item $V = Aw,$ 
\item $X_iw = \zeta_i w$ for  all $i \in \mathcal I$.
\end{enumerate}
We refer to the pair $(V,w)$ as a {\it Whittaker pair of type $\zeta$}.  

In what follows,  any  $v \in V$ such that $X_iv = \zeta_i v$ for all $i \in \mathcal I$ will be called a {\it Whittaker vector of type $\zeta$}.   Such a vector is simply a common eigenvector for all the generators $X_i$ with nonzero eigenvalues.    For a Whittaker module $V$ with cyclic Whittaker vector $w$ of type $\zeta$, let $Q = \ann_R (w)$, the annihilator of $w$ in $R$.  Note that $Q$ is a left ideal of $R$ and  $\ann_A (w)$ is a left ideal of $A$,   while $\ann_A (V)$ is an ideal of $A$.   We fix this notation for the remainder of the paper. 

In Section \ref{sec:univeralObject} we construct a universal Whittaker module of type $\zeta$ for each generalized Weyl algebra $A = R(\phi,t)$.   This module is used in the proof of Theorem \ref{thm:bijectIdealsModules} to show that the isomorphism classes of Whittaker modules of type $\zeta$ are in bijection with the $\phi$-stable left ideals of $R$.  In particular, simple Whittaker modules correspond to maximal $\phi$-stable left ideals of $R$.   For finite-dimensional complex semisimple Lie algebras $\mathfrak g$,  the corresponding result  in \cite{{kostant:wvrt78}}  states that the isomorphism
classes of Whittaker modules of type $\zeta$ are in bijection with the
ideals of the center of the universal enveloping algebra $U(\mathfrak g)$.  A similar result holds for the quantum enveloping algebra $U_q(\mathfrak {sl}_2)$ (see  \cite{{ondrus:thesis}, {ondrus:Whittaker}}).      For an arbitrary Whittaker module $V = Aw$ for a generalized Weyl algebra $A$,  in Section 4  we obtain a description of  the annihilator of $w$:   \    $\ann_A(w)  =  AQ  + \sum_{i \in \mathcal I}  A(X_i - \zeta_i)$, where $Q = \ann_R(w)$.  In Sections 5 and 6,  we impose the assumption  that $R$ is commutative and determine the Whittaker vectors inside a Whittaker module.    When $R$ is commutative, $Q$ is a prime ideal not containing any $t_i$, and the center of $A$ is contained in $R$,   then   $\ann_A(V) = AQ$ by Theorem  \ref{thm:annVwhenQprime}.    The final sections are devoted to illustrating what these results say for certain well-known algebras such as the (quantum) Weyl algebra, the quantum plane, and Smith's generalizations of $U(\mathfrak{sl}_2)$ and of $U_q(\mathfrak{sl}_2)$.  We recover the results of  \cite{{tang:WMSmithAlg}} and \cite{JiWangZhou:quantumSmith} for the (quantum) Smith algebras of characteristic zero and determine the Whittaker modules for all these algebras in the  modular and root of unity cases.

\begin{center} {\textbf{Acknowledgments}}   \end{center}
This paper was completed while the first author was a Simons Visiting Professor 
in the Combinatorial Representation Theory Program at 
the Mathematical Sciences Research
Institute.  She  acknowledges with gratitude  the warm  hospitality of MSRI and  the support from 
 the Simons Foundation,   as well as 
  from 
National Science Foundation grant \#{}DMS--0245082.  
The authors also  thank Louis Solomon and Samuel Lopes  for their insights and helpful comments.

\section{Basic facts about generalized Weyl algebras \\
and their Whittaker modules} \label{sec:basicFacts}

Assume  $A = R(\phi, t)$ is a generalized Weyl algebra as in Section  \ref{sec:notationGWA}.  Let $\Gamma$ denote the semigroup of tuples $\gamma = (\gamma_i)_{i \in \mathcal I}$  of nonnegative integers with only finitely many nonzero entries   under componentwise addition,  $\gamma + \delta = \big((\gamma+\delta)_i\big)_{i \in\mathcal I}$ where $(\gamma+ \delta)_i = \gamma_i + \delta_i$.   For $\gamma \in \Gamma$, set 
\begin{equation}\label{eqn:XtoPower}
X^\gamma =  \prod_{i \in \mathcal I}  X_i^{\gamma_i}. 
\end{equation}
Because the various $X_j$ commute,  it follows that  $X^\gamma X^\delta = X^{\gamma + \delta}$.

We adopt the notation

$$Z_i^\ell= \left\{ \begin{array}{ll} X_i^\ell & \mbox{if $\ell \ge 0$} \\
Y_i^{-\ell} & \mbox{if $\ell < 0$.}  \end{array} \right.$$
Observe that  the defining relations give

\begin{equation}\label{eq:express}  Y_i ^k X_i^\ell   
= \left\{ \begin{array}{ll} \phi_i^{-(k-1)}(t_i) \phi_i^{-(k-2)}(t_i) \cdots \phi_i^{-(k-\ell)}(t_i)Y_i^{k-\ell} & \mbox{if $k \ge \ell$} \\
\phi_i^{-(k-1)}(t_i) \phi_i^{-(k-2)}(t_i) \cdots \phi_i^{-1}(t_i) \,t_iX_i^{\ell-k}& \mbox{if $k < \ell$.}  \end{array} \right. \end{equation}  
View $\mathbb Z^{\mathcal I}$ under componentwise addition, and  let $\Lambda$ denote the subgroup of  $\mathbb Z^{\mathcal I}$  of all tuples $\alpha = (\alpha_i)_{i \in \mathcal I}$ having only finitely many nonzero components.   Set 
$$Z^\alpha =  \prod_{i}  Z_i^{\alpha_i},$$ 
and note that this product is well-defined since the $Z_i$ commute.  Then it follows from \eqref{eq:express}  that every element  $a \in A = R(\phi, t)$ can be written as a finite sum

$$a = \sum_{\alpha \in \Lambda}  c_\alpha Z^{\alpha}$$
with coefficients $c_\alpha$ in  $R$.

\begin{lem}\label{lem:Rfree}  $A = R(\phi, t)$ is a free left (or right) $R$-module with 
basis $\{Z^\alpha \mid \alpha \in \Lambda\}$.   \end{lem}

\begin{proof}  We have observed already  that these elements span $A$ over $R$.
Now suppose that  $a = \sum_{\alpha}  c_\alpha Z^{\alpha} = 0$, where $c_\alpha \in R$ for
all $\alpha \in \Lambda$.    Given such an expression,  for each $i \in\mathcal I$ set  

$$\gamma_i = \max\Big(\{ -\alpha_i  \mid \alpha_i < 0 \} \cup \{0\}\Big),$$

\noindent where the maximum is taken over all $\alpha$ such that $c_\alpha \neq 0$.
Then  $\gamma = (\gamma_i) \in \Gamma$, and  by \eqref{eq:express} we have

\begin{equation}\label{eq:Xpoly}  a X^\gamma = \sum_{\alpha \in \Lambda} c_\alpha d_\alpha Z^{\alpha + \gamma}   \end{equation}
 for some nonzero $d_\alpha \in R$.   The powers occurring in the monomials
 $Z^{\alpha+\gamma}$ are all nonnegative.   Thus, the factors in  $Z^{\alpha+\gamma}$ are just 
 the $X_i$.   
 Because the subalgebra of $A$ generated by $R$ and the $X_i$ is a skew-polynomial
 ring, it is free over $R$ with basis the monomials in the $X_i$.   Therefore,  
  $c_\alpha d_\alpha = 0$, and hence   $c_\alpha = 0$   for all $\alpha \in \Lambda$.  \end{proof}

The next proposition  is a generalization of a result of Kulkarni \cite[Cor.~2.02]{kulkarni:jalg01} which
treats the case that $| {\mathcal I} | = 1$.

\begin{prop}\label{prop:Acent}  Let $A = R(\phi,t)$ be a generalized Weyl algebra with $R$ commutative.
Then the center $ \mathcal Z(A)$ of $A$  is generated by
the
elements of $R$ in  $R^\phi:=  \{r \in R \mid \phi_i(r) = r$ for all $i \in \mathcal I\}$
and all  the monomials $Z^\alpha$ for   $\alpha \in \Lambda$ 
such that $\phi^\alpha: = \prod_{i \in \mathcal I} \phi_i^{\alpha_i} = \hbox{\rm id}_R$.    \end{prop}

\begin{proof}  It is easy to see from (\ref {eqn:Xir}) that  $R^\phi$
is contained in $\mathcal Z = \mathcal Z(A)$.    Moreover,  the relations in \eqref{eq:express}  can
be used to show that $Z^\alpha \in \mathcal Z$ if and only if $\phi^\alpha = \hbox{\rm id}_R$.  
Now suppose $\sum_{\alpha \in \Lambda} r_\alpha Z^\alpha  \in \mathcal Z$.   Then for $s \in R$, 
\begin{equation} \sum_{\alpha \in \Lambda} s r_\alpha Z^\alpha =   s\left(\sum_{\alpha \in \Lambda} r_\alpha Z^\alpha\right)  =  \left(\sum_{\alpha \in \Lambda} r_\alpha Z^\alpha\right) s  
= \sum_{\alpha \in \Lambda} r_\alpha \phi^\alpha(s) Z^\alpha.
\end{equation}
Thus if $r_\alpha \neq 0$, then $s = \phi^\alpha(s)$ for all $s \in R$ by
Lemma \ref{lem:Rfree},   so that $\phi^\alpha = \hbox{\rm id}_R$ and  $Z^\alpha \in \mathcal Z$. 
But then $ \left ( \sum_\alpha r_\alpha Z^\alpha \right) X_i  =
X_i \left ( \sum_\alpha r_\alpha Z^\alpha \right)$ implies 
$\sum_\alpha r_\alpha X_i Z^\alpha = \sum_\alpha \phi_i(r_\alpha) X_i Z^\alpha$ since
$Z^\alpha \in \mathcal Z$ for each nonzero $r_\alpha$.    This forces
$r_\alpha = \phi_i(r_\alpha)$ for all $i$, so that  $r_\alpha \in R^\phi$.  
  \end{proof}

\begin{lem}\label{lem:V=Rw}
Let $V$ be a Whittaker module for $A$ with cyclic Whittaker vector $w$ of type $\zeta$.  Then $V = Rw$.  If $R$ is commutative, then $\ann_R(V) = \ann_R(w)$.
\end{lem}
\begin{proof}
Observe that $X_iw = \zeta_i w \in Rw$ and $Y_iw = \zeta_i^{-1} Y_iX_iw = \zeta_i^{-1} t_iw \in Rw,$ and then apply the relations $X_ir = \phi_i (r) X_i$ and $Y_ir = \phi_i^{-1}(r)Y_i$ for $r \in R.$  The assertion about annihilators is an immediate consequence of the fact that $V = Rw$ and the commutativity
of $R$. 
\end{proof}

\begin{defn}
If $J \subseteq R$ is an ideal or left ideal of $R,$ we say that $J$ is \emph{$\phi$-stable} if $\phi_i (J) \subseteq J$ for all $i \in \mathcal I$. 
\end{defn}

\begin{exams} For a fixed Whittaker module $V = Aw$, if $r \in Q = \ann_R (w),$ then $0 = X_irw = \phi_i (r) X_iw = \zeta_i \phi_i (r)w,$ so it follows that $Q$ is a $\phi$-stable left ideal of $R.$  For another example, assume $z \in R$ is fixed by $\phi_i$ for all $i$.   Then the left ideal $(z) = Rz$ of $R$ generated by $z$ is clearly $\phi$-stable.   \end{exams}

\begin{rem}\label{rem:radicalStable}
The sum of two $\phi$-stable (left) ideals is again $\phi$-stable.  In addition, if $R$ is commutative and $J$ is a $\phi$-stable ideal of $R$, it is easily seen that the radical $\sqrt J$ is also a $\phi$-stable ideal of $R$. 
\end{rem}

If $J$ is a $\phi$-stable left ideal of $R$,  applying $\phi_i^{-1}$ to the containment $\phi_i (J) \subseteq J$ gives  $J \subseteq \phi_i^{-1}(J)$.   Repeating this indefinitely yields $J \subseteq \phi_i^{-1}(J) \subseteq \phi_i^{-2}(J) \subseteq \phi_i^{-3}(J) \subseteq \cdots$.    Since the $\phi_i^{-k}(J)$ ($k \ge 0$) form an ascending chain of left ideals of the Noetherian ring $R$,  it follows that $\phi_i^{-k}(J) = \phi_i^{-(k+1)}(J)$ for some $k \ge 0.$  Applying an appropriate power of $\phi_i$ to each side gives the following lemma.

\begin{lem}\label{lem:phiStable}
If $J$ is a $\phi$-stable left ideal of $R,$ then $\phi_i (J) = J = \phi_i^{-1}(J)$ for all $i.$
\end{lem}

Let $J$ be a $\phi$-stable left ideal of $R$ such that $t_i  \not\in J$ for all $i \in \mathcal I$.  By $\overline r$, we mean the coset $\overline r = r+J$ for all $r \in R$.  Thus, $\overline {t_i}  \neq 0$ for $i \in \mathcal I$.  Since $J$ is $\phi$-stable, we have an induced map $\overline {\phi_i}$ on the quotient $R/J$ given by $\overline {\phi_i}(\overline r) = \phi_i(r) + J = \overline{\phi_i(r)}$.  When  $R$ is commutative, then $R/J$ is a ring, and $\overline {\phi_i}$ is an automorphism of $R/J$ by Lemma \ref{lem:phiStable}. Thus the following result holds in that case.

\begin{prop} \label{prop:quotientGWA}
Let $A = R(\phi, t)$ be a generalized Weyl algebra.  Assume $R$ is commutative, and let $J$ be a $\phi$-stable prime ideal of $R$ such that $t_i \not \in J$ for all $i \in \mathcal I$.    Then  $A/AJ$ is isomorphic to the generalized Weyl algebra  $\overline A := (R/J)(\overline \phi, \overline t)$.   \end{prop}

\begin{proof}  We will write bars on the $X_i$ and $Y_i$ in $\overline A$ to distinguish them from the generators in $A$, although the bar does not denote a coset reduction in this instance.  Consider the $\kk$-algebra homomorphism  $\Phi: { \mathcal F}  \to \overline A$  from the free algebra $\mathcal F$ generated by $R$, $X_i, Y_i$, $i \in \mathcal I$, to $\overline A$  given  $r \mapsto \overline r$, $X_i \mapsto \overline {X}_i$,  and $Y_i \mapsto \overline{Y}_i$.      Then $\Phi(Y_i X_i) = \overline{Y}_i\,{\overline X}_i = \overline{t_i} = \Phi(t_i)$, so that $Y_i X_i - t_i \in \ker \Phi$.     Similarly,  $\Phi(X_i r) = \overline{X}_i\overline {r} = \overline {\phi_i}({\overline r}) \overline{X}_i =   \overline{\phi_i(r)} \overline{X}_i  = \Phi(\phi_i(r) X_i)$, and $X_i r - \phi_i(r)X_i \in \ker \Phi$.  Arguing in this way, we see that  there is an induced algebra homomorphism $\overline \Phi: A \to \overline A$.   Clearly, $AJ$ is in the kernel.  Now if some $\sum_{\alpha \in \Lambda}  Z^\alpha  r_\alpha$  maps to $0$, then $\sum_{\alpha \in \Lambda}\overline{Z}^\alpha \overline{r_\alpha}  = 0$ in $\overline A$, so by the freeness of $\overline A$ as a module for the domain $R/J$  (see Lemma \ref{lem:Rfree}), we obtain $\overline {r_\alpha}= 0$ for each $\alpha \in \Lambda$.   This implies that $r_\alpha  \in J$ for all $\alpha$, so that $\sum_{\alpha \in \Lambda}  Z^\alpha r_\alpha  \in AJ$.  \end{proof}
 
\section{Constructing a universal object} \label{sec:univeralObject}

We continue to assume that $A = R(\phi,t)$ is a generalized Weyl algebra.

\begin{defn}\label{def:univWhitt}
Let $(V,w)$ be a Whittaker pair of type $\zeta$.  Suppose that $V$ has the property that for any other Whittaker pair $(V',w')$ of type $\zeta$, there exists a unique surjective module homomorphism 
$\sigma' : V \to V'$ such that $\sigma' (w) = w'$.  Then we say that $(V,w)$ is a \emph{universal Whittaker pair of type $\zeta$} and  $V$ is a \emph{universal Whittaker module of type $\zeta$}.
\end{defn}

Suppose that $(V_1,w_1)$ and $(V_2,w_2)$ are universal Whittaker pairs of type $\zeta$.  Then there are surjective $A$-module homomorphisms $\sigma_2 : V_1 \to V_2$ and $\sigma_1 : V_2 \to V_1$ such that $\sigma_2(w_1) = w_2$ and $\sigma_1(w_2) = w_1$.  If $v \in V_1$, then we may write $v = rw_1$, and thus $\sigma_2(v) = \sigma_2(rw_1) = r \sigma_2(w_1) = rw_2$.  Moreover, $\sigma_1(rw_2) = r \sigma_1(w_2) = rw_1 = v$, so we see that $\sigma_2 \circ \sigma_1 = {\rm id}_{V_1}$.  Similarly,  $\sigma_1 \circ \sigma_2 = {\rm id}_{V_2}$.  Thus the maps $\sigma_1$ and $\sigma_2$ are isomorphisms of $A$-modules, and it makes sense to refer to a universal Whittaker module of type $\zeta$ as \emph{the} universal Whittaker module of type $\zeta$.

To construct a universal Whittaker pair $(V_{\mathfrak u},w_{\mathfrak u})$ of type $\zeta$, we define an action of $A$ on $R$ via 
\begin{align}
r'.r  &= r'r, \label{eqn:univModueRaction} \\
X_i.r  &= \zeta_i \phi_i(r), \quad \mbox{and}  \label{eqn:univModueXaction} \\
Y_i.r  &= \zeta_i^{-1}\phi_i^{-1}(r)t_i \label{eqn:univModueYaction}
\end{align}
for $r, r' \in R$ and $i \in \cal I$.  It is straightforward to verify that under this action, $R$ is a Whittaker module of type $\zeta$ with cyclic Whittaker vector  $1$.  When we regard $R$ as a Whittaker module with the above action, we write $V_{\mathfrak u} = R$ and $w_{\mathfrak u} = 1$.

\begin{lem} \label{lem:universalWhitt}
The module $V_{\mathfrak u}$ is the universal Whittaker module of type $\zeta$  and   $\ann_R (w_{\mathfrak u}) =0$.
\end{lem}
\begin{proof} It is clear that $\ann_R(w_{\mathfrak u}) = 0$.  Let $(V,w)$ be an arbitrary Whittaker pair of type $\zeta$,  and define a map $\sigma : V_{\mathfrak u} \to V$ as follows.  For $v \in V_{\mathfrak u}$, set $\sigma (v) = rw \in V$, where $r \in R$ is such that $v = rw_{\mathfrak u}$.  If $s \in R$ satisfies $sw_{\mathfrak u} = v = rw_{\mathfrak u}$, then $s-r \in \ann_R (w_{\mathfrak u}) = 0$, so $s = r$ and the map $\sigma : V_{\mathfrak u} \to V$ is well-defined.

With $v = rw_{\mathfrak u} \in V_{\mathfrak u}$, we must verify that $\sigma (av) = a \sigma (v)$ for all $a \in A$.  But since $A$ is generated over $R$ by $X_i$ and $Y_i$ for $i \in \mathcal I$,  it is sufficient to consider the cases $a \in R$, $a = X_i$, and $a = Y_i$, and these routine calculations are omitted.

Because $V = Aw = Rw$, we have  that $\sigma (V_{\mathfrak u}) = \sigma (Rw_{\mathfrak u}) = R \sigma(w_{\mathfrak u}) = Rw = V$, and thus $\sigma : V_{\mathfrak u} \to V$ is surjective.  The uniqueness of $\sigma$ follows from the fact that $\sigma (w_{\mathfrak u}) = w$ and $\sigma$ respects the action of $R$.
\end{proof}

\begin{rem}  An $\kk$-basis for $V_{\mathfrak u}$ is the set  $\{b_\ell w_{\mathfrak u}\}_{\ell \in \mathcal L}$, where $\{b_\ell\}_{\ell \in \mathcal L}$ is any  $\kk$-basis of $R$.   \end{rem}

\begin{rem}  Here, we describe an alternative construction of the universal Whittaker module of type $\zeta$ similar to that of \cite[Thm.~3.3]{kostant:wvrt78}.  As a convenient shorthand in the construction,  let $\kk[X]$ denote the polynomial algebra over $\kk$ generated by the $X_i$, $i \in \mathcal I$, and regard  $\kk [X]$ as a subalgebra of $A$.     Give the one-dimensional space $\kk w_\zeta$ an $\kk [X]$-module structure according to action $X_i w_\zeta = \zeta_i w_\zeta$.  Set 
\begin{equation}\label{eqn:constructUniv}
V_\zeta = A \otimes_{\kk[X]} \kk w_\zeta,
\end{equation}
and (making a slight abuse of notation) write $w_\zeta$ to denote $1 \otimes w_\zeta$.  Then it is clear that $V_\zeta = Aw_\zeta$ and $X_iw_\zeta = \zeta_iw_\zeta$, so $(V_\zeta ,w_\zeta)$ is a Whittaker pair of type $\zeta$.    That this induced construction also gives the universal Whittaker module
follows from the fact that the subalgebra of $A$  generated by $R$ and the $X_i$ is
free over $\kk[X]$.  We omit the details.  
\end{rem}

\begin{lem} \label{lem:submodulesOfR}
The $A$-submodules of $V_{\mathfrak u} = R$ are exactly the $\phi$-stable left ideals of $R$.
\end{lem}
\begin{proof}
Suppose that $J \subseteq R$ is a submodule of $V_{\mathfrak u}$.  Equation (\ref{eqn:univModueRaction}) shows  that $J$ is a left ideal of $R$.  Since $\zeta_i$ is nonzero for all $i \in \cal I$,  (\ref{eqn:univModueXaction}) implies that $\phi_i(J) \subseteq J$ for all $i$, and thus $J$ is $\phi$-stable.   It is routine to verify that any $\phi$-stable left ideal of $R$ is a submodule of $V_{\mathfrak u}$.   
\end{proof}

\begin{defn} \label{defn:V_Q}
If $Q$ is a $\phi$-stable ideal of $R$, let $V_Q = R/Q$, and regard $V_Q$ as the quotient $V_{\mathfrak u} / QV_{\mathfrak u}$ with cyclic Whittaker vector $w_Q = 1+Q$.
Observe that $\ann_R (w_Q) = \{ r \in R \mid r(1+Q) = 0+Q \} = Q$.  
\end{defn}

Suppose now  that $(V,w)$ is an arbitrary Whittaker pair of type $\zeta$, and let $Q = \ann_R (w)$. 
Then there is a map $\sigma : V_{\mathfrak u} \to V$,  $r w_{\mathfrak u} \mapsto r w$.
 If $0 = \sigma (r w_{\mathfrak u}) = rw$, then  $r \in Q$, and thus $v = rw_{\mathfrak u} \in Qw_{\mathfrak u}$.  Hence  $\ker (\sigma) = Qw_{\mathfrak u}$ and $V \cong V_{\mathfrak u} / Qw_{\mathfrak u} = V_Q$.  We therefore have

\begin{lem} \label{lem:Whittquo} 
Assume $(V,w)$ is an arbitrary Whittaker pair of type $\zeta$, and let $Q = \ann_R (w)$.  Then $V \cong V_{\mathfrak u} / Qw_{\mathfrak u} = R/Q = V_Q$, where $(V_{\mathfrak u}, w_{\mathfrak u})$ is the universal Whittaker pair of type $\zeta$.   
\end{lem}
 
Theorem \ref{thm:bijectIdealsModules} below  is the generalized Weyl algebra analogue  of Kostant's result  \cite [Thm.~3.2]{ kostant:wvrt78}  for finite-dimensional complex semisimple Lie algebras  and Ondrus' result \cite[Cor.~4.1]{ondrus:Whittaker} for the quantum group $U_q(\mathfrak{sl}_2)$.

\begin{thm}\label{thm:bijectIdealsModules}   Let $A = R(\phi,t)$ be a generalized Weyl algebra.   Then the map 
$$ \left\{ {\mbox{isomorphism classes of} \atop \mbox{Whittaker pairs of type $\zeta$}} \right\}    \psito \{ \mbox{$\phi$-stable left ideals of $R$} \}$$
given by 
$$(V,w) \mapsto \ann_R(w)$$
is a bijection.
\end{thm}
\begin{proof}   Suppose that $(V_1, w_1)$ and $(V_2, w_2)$ are Whittaker pairs of type $\zeta$ with $\ann_R(w_1) = \ann_R(w_2)$, and set $Q = \ann_R(w_1) = \ann_R(w_2)$.  In Lemma \ref{lem:Whittquo}, we have seen that $V_j \cong V_{\mathfrak u} / Q w_{\mathfrak u}$ for $j = 1,2$, where $(V_{\mathfrak u},w_{\mathfrak u})$ is the universal Whittaker pair of type $\zeta$.   Thus $V_1 \cong V_2$.   This implies that the map $\Psi: (V,w) \mapsto \ann_R(w)$ is injective.

Now suppose $Q \subseteq R$ is a $\phi$-stable left ideal of $R$, and let $V_Q$ be as in Definition \ref{defn:V_Q}.  Since  $\ann_R (w_Q) = Q$,  the map   $\Psi$ is surjective as well.   \end{proof}

\begin{cor} \label{cor:whenIsVsimple}
Let $(V,w)$ be a Whittaker pair of type $\zeta$ for a generalized Weyl algebra $A = R(\phi,t)$.  Then $V$ is simple if and only if $\ann_R(w)$ is a maximal $\phi$-stable left ideal of $R$.
\end{cor}

\begin{exam} 
Fix an element $q \in \kk$ with $q \neq 0$ and $q^2 \neq 1$.  Let  $R = \kk [c,K,K^{-1}]$ and $t = c-\frac{qK+q^{-1}K^{-1}}{(q-q^{-1})^2} \in R$, and define $\phi : R \to R$ by $K \mapsto q^{-2}K$ and $c \mapsto c$.  (Because  $ |{\mathcal I} | = 1$ in this example, we are omitting the subscripts on $\phi$ and $t$.)  Then $A = R(\phi,t) \cong U_q(\mathfrak{sl}_2)$.  Since $R$ is commutative, the simple Whittaker modules correspond to maximal $\phi$-stable (two-sided) ideals of $R$.  If $\xi \in \kk$, then the ideal $R(c - \xi)$ generated by $c - \xi$ is clearly $\phi$-stable.  We shall see in Section \ref{sec:qSmith} that this is a maximal $\phi$-stable ideal when $q^2$ is not a root of unity.
\end{exam}

\begin{exam}\label{ex:nthWeyl}
Assume  $\kk$ has characteristic 0, and let $A_n = R_n(\phi,t)$ where $R_n = \kk[t_1,\dots,t_n]$ and $\phi_i (t_j) = t_j - \delta_{i,j}$.   Assume $Y_iX_i = t_i$, $X_iY_i = \phi_i(t_i)$, $[X_i,X_j] = 0 = [Y_i,Y_j]$,  and $[X_i,Y_j] = 0$ for $i \neq j$.   Then $A_n$ is the $n$th Weyl algebra realized as a generalized Weyl algebra.  It is straightforward to show that $R_n$ contains no  proper $\phi$-stable ideals, and thus every Whittaker module $V$ for $A_n$ is simple.  In particular, the universal Whittaker module $V_{\mathfrak u}$ of type $\zeta$  is simple and is the unique Whittaker module of type $\zeta$ for  $A_n$.  The set  $\{ t^\gamma = t^\gamma w_{\mathfrak u} \mid  t^\gamma = \prod_{i =1}^n t_i^{\gamma_i},  \, \gamma_i \in \mathbb Z_{\geq 0} \}$ is a basis for $V_{\mathfrak u}$,
and the $A_n$-action on $V_{\mathfrak u}$ is given by
\begin{eqnarray}\label{eq:Weylex}  t^\beta.t^\gamma &=&  t^{\beta+\gamma}, \nonumber \\
X_i.t^\gamma &=&  \zeta_i (t_i-1)^{\gamma_i}\prod_{j \neq i} t_j^{\gamma_j} \\
Y_i.t^\gamma &=&   \zeta_i^{-1} t_i(t_i+1)^{\gamma_i}\prod_{j \neq i} t_j^{\gamma_j}. \nonumber \end{eqnarray}
\end{exam}

If  $\mathcal K$ is an ideal of the center $\mathcal Z = \mathcal Z(A)$ of a generalized Weyl algebra $A = R(\phi,t)$, then $\mathcal K V_{\mathfrak u}$ is a submodule of the universal Whittaker module   $V_{\mathfrak u} =   R$.   Our next goal is to show that when $R$ is commutative and  $\mathcal I$ is finite,  then under some assumptions, 
\begin{equation}\label{eq:VuK} V_{\mathfrak u, \mathcal K} : = V_{\mathfrak u} \Big /  \mathcal K V_{\mathfrak u}  \end{equation} is simple for every maximal ideal $\mathcal K$ of $\mathcal Z$. 

Recall that a commutative ring is said to be a \ {\it Jacobson ring} \ if each prime ideal is
the intersection of maximal ideals.      We will use the following two results:

\begin{thm}{\rm \cite[Thm.~4.19]{eisenbud:CommAlg}} \label{thm:JacRing}
Let $S$ be a Jacobson ring.  If $T$ is a finitely generated $S$-algebra, then $T$ is a Jacobson ring.  Furthermore, if $\mathfrak n \subseteq T$ is a maximal ideal, then $\mathfrak m := \mathfrak n \cap S$ is a maximal ideal of $S$, and $T / \mathfrak n$ is a finite extension field of $S / \mathfrak m$.
\end{thm}

\begin{thm}{\rm \cite[Thm.~6.20]{swan:KtheoryFiniteGroups}}\label{thm:finGenMonoid}
Let M be a finitely generated commutative monoid and let $\{ f_i:  M \to \z  \mid i = 1, ... ,n \}$ be a finite collection of homomorphisms.  Then $G = \{ x \in M \mid  \mbox{$f_i(x) \geq 0$ for all $i$} \}$ is a finitely generated monoid.
\end{thm}

Any field $\kk$ is a Jacobson ring, and hence by Theorem \ref{thm:JacRing},  so is any
finitely generated commutative $\kk$-algebra.  We intend to apply this theorem to the pair
$S = R^\phi$ and $T = \mathcal Z = \mathcal Z(A)$, where our notation is that of Proposition \ref{prop:Acent}.     Thus,  we need conditions under which $\mathcal Z$ is a finitely generated $R^\phi$-algebra.  

Let $\Delta = \{ \alpha \in \Lambda \mid \phi^\alpha = {\rm id}_R \}$, and note that $\mathcal Z = \bigoplus_{\alpha \in \Delta} R^\phi Z^\alpha$ by Proposition  \ref{prop:Acent}.  If $|\mathcal I|$ is finite, then the subgroup $\Delta \subseteq \Lambda$ is finitely generated.  However, it may not be the case that $Z^\alpha Z^\beta = Z^{\alpha + \beta}$, so it is not immediately obvious that $\mathcal Z$ is a finitely generated $R^\phi$-algebra.  \medskip

\begin{lem} \label{lem:ZfiniteOverRphi}
If $R$ is commutative and $| {\mathcal I} | < \infty$, then $\mathcal Z$ is a finitely generated $R^\phi$-algebra. 
\end{lem}
\begin{proof}
Assume that $|{\mathcal I}| = n < \infty$ so that $\Lambda = \z^n$, and let $\Sigma = \{ \pm 1 \}^n$.  For $\varepsilon = ( \varepsilon_1, \ldots, \varepsilon_n) \in \Sigma$, define homomorphisms $f_i^\varepsilon: \z^n \to \z$ (for $i = 1, \ldots, n$) by $f_i^\varepsilon ( \alpha )= \varepsilon_i \alpha_i$.  With $\Delta = \{ \alpha \in \Lambda \mid \phi^\alpha = {\rm id}_R \} \subseteq \z^n$ as above, let $\Delta_\varepsilon \subseteq \Delta$ be the monoid defined by $\Delta_\varepsilon = \{ \alpha \in \Delta \mid \mbox{$f_i^\varepsilon (\alpha) \ge 0$ for $i = 1, \ldots, n$} \}$.  Note that $Z^\alpha Z^\beta = Z^{\alpha + \beta}$ for $\alpha, \beta \in \Delta_\varepsilon$, and by Theorem \ref{thm:finGenMonoid}, there is a finite set $\mathcal G_\varepsilon$ of generators for the monoid $\Delta_\varepsilon$.  Observe that $|\Sigma| = 2^n$, and $\Delta = \bigcup_{\varepsilon \in \Sigma} \Delta_\varepsilon$.  Thus the set $\mathcal G = \bigcup_{\varepsilon \in \Sigma} \mathcal G_\varepsilon$ is finite, and the set $\{ Z^\alpha \mid \alpha \in \mathcal G \}$ is a finite set of generators for $\mathcal Z$ over $R^\phi$.
\end{proof}

\begin{defn}
We say that a $\phi$-stable ideal $Q$ of $R$ is {\it centrally generated} if $Q = R (Q \cap \mathcal Z)$.  
\end{defn}

\begin{lem}\label{lem:QintZ}
Let $A$, $\mathcal Z$, $\mathcal K$, and $V_{\mathfrak u, \mathcal K}$ be as above and set 
 $$w_{\mathfrak u,\mathcal K} = 1 + \mathcal K V_{\mathfrak u} \in V_{\mathfrak u,\mathcal K} 
 = V_{\mathfrak u}/\mathcal K V_{\mathfrak u}.$$  
 \noindent Assume $\mathcal K$ is a maximal ideal of $\mathcal Z$, and let $Q = \ann_R(w_{\mathfrak u, \mathcal K})$.  Then $Q \cap \mathcal Z = R^\phi \cap \mathcal K$,
 where $R^\phi = R \cap \mathcal Z = \{ r \in R \mid \mbox{$\phi_i(r) = r$ for all $i \in \mathcal I$} \}$. 
\end{lem}

\begin{proof}  It follows from  the construction of $V_{\mathfrak u, \mathcal K}$ that  $\mathcal K \subseteq \ann_{\mathcal Z}(w_{\mathfrak u, \mathcal K})$.  However,  $\ann_{\mathcal Z}(w_{\mathfrak u, \mathcal K})$ is clearly a proper ideal of $\mathcal Z$,  so since $\mathcal K$ is maximal,  
$\mathcal K = \ann_{\mathcal Z}(w_{\mathfrak u, \mathcal K})$ must hold.  The proof of the remaining assertions  is straightforward.    \end{proof}

\begin{thm}\label{thm:Whittsimple}
Assume $R$ is commutative and every maximal $\phi$-stable ideal of $R$ is centrally generated.  Let $\mathcal K$ be a maximal ideal of  the center $\mathcal Z$ of $A= R(\phi,t)$.  If $|{\mathcal I}| < \infty$ and $R^\phi$ is a finitely generated $\kk$-algebra, then the Whittaker module $V_{\mathfrak u, \mathcal K}$ is simple.
Moreover,  if $Q = \ann_R(w_{\mathfrak u, \mathcal K})$, where $w_{\mathfrak u, \mathcal K} =
1+\mathcal K V_{\mathfrak u}$,   then $Q = R(R^\phi \cap \mathcal K)$.

\end{thm}
\begin{proof}
Since $R^\phi$ is a finitely generated $\kk$-algebra, it  is a Jacobson ring.  Consequently,  by 
Lemma  \ref{lem:ZfiniteOverRphi} and Theorem \ref{thm:JacRing},  $\mathcal Z$ is a
finitely generated $R^\phi$-algebra, and $R^\phi \cap \mathcal K$ is a maximal ideal of $R^\phi$.  Let $Q = \ann_R(w_{\mathfrak u, \mathcal K})$, and recall from Lemma \ref{lem:QintZ} that $Q \cap \mathcal Z = R^\phi \cap \mathcal K$.   Since $R$ is Noetherian, there exists a maximal $\phi$-stable ideal  $Q'$ of $R$ containing $Q$.  By assumption $Q' = R(Q' \cap \mathcal Z)$.  But $Q' \cap \mathcal Z$ is a proper ideal of $R^\phi$ because $1 \not\in Q' \cap \mathcal Z$, and $Q' \cap \mathcal Z \supseteq Q \cap \mathcal Z$.  As $Q \cap \mathcal Z = R^\phi \cap \mathcal K$ is a maximal ideal of
$R^\phi$, it follows that $Q' \cap \mathcal Z = Q \cap \mathcal Z$, and so  $$Q' = R(Q' \cap \mathcal Z) = R(Q \cap \mathcal Z) \subseteq  Q.$$
This implies that $Q = Q'$ is maximal among $\phi$-stable ideals of $R$, hence $V_{\mathfrak u, \mathcal K}$ is simple by Corollary \ref{cor:whenIsVsimple}.   But then $Q = R(Q \cap \mathcal Z) = R (R^\phi \cap \mathcal K)$, as claimed.  
\end{proof}

\begin{rem} {All the examples in Sections 8-10 satisfy the hypothesis that  $R$ is commutative and $R^\phi$ is a finitely generated $\kk$-algebra.   Many  of the examples satisfy the condition that every maximal $\phi$-stable ideal of $R$ is centrally generated, and thus the module $V_{\mathfrak u, \mathcal K} = V_{\mathfrak u}/\mathcal K V_{\mathfrak u} \cong V_{\mathfrak u}/Q V_{\mathfrak u} = R/Q$, where
$Q = \ann_R(w_{\mathfrak u,\mathcal K})$,   is a simple Whittaker module in those cases.}  \end{rem}

\section{An expression for $\ann_A(w)$}

Let $A = R(\phi,t)$ be a generalized Weyl algebra as in 
Section \ref{sec:notationGWA},   and suppose that $V = Aw$ is a Whittaker module of type $\zeta$ with $Q = \ann_R(w)$.  The map $A \to V$ given by $a \mapsto aw$ shows that $V \cong A/\ann_A(w)$, and it is clear that $AQ + \sum_{i \in \mathcal I} A(X_i-\zeta_i) \subseteq \ann_A(w)$.  In this section, we prove that in fact these two left ideals of $A$ always coincide.

As before, let  $\Gamma$ denote the semigroup of  
tuples $\gamma = (\gamma_i)_{i \in \mathcal I}$  of nonnegative integers with only
finitely many nonzero entries  under componentwise addition, and let $X^\gamma =  \prod_{i \in \mathcal I}  X_i^{\gamma_i}$.  
For  $\gamma \in \Gamma$,  set  $| \gamma | = \sum_{i \in \mathcal I} \gamma_i$.

\begin{lem}\label{lem:aXm}
Let $I = AQ + \sum_{i\in \mathcal I} A(X_i-\zeta_i)$, where $(V,w)$ is a Whittaker pair of type $\zeta$ and $Q = \ann_R(w)$.  Let $a \in A$, and suppose that there exists $\gamma \in \Gamma$ such that $aX^\gamma \in I$.  Then $a \in I$.
\end{lem}
\begin{proof}
The proof is by induction on $|\gamma|$.   We may assume that $|\gamma| > 0$ since there is nothing to prove if $\gamma_i = 0$ for all $i$.   Assume that $\gamma_k>0$ for some $k \in \mathcal I$, and thus we use the assumption that $aX^\gamma \in I$ to show that $aX^{\gamma'} \in I$, where $\gamma' = (\gamma'_i)$  is such that  $\gamma'_i = \gamma_i$ for $i \neq k$ and
$\gamma'_k = \gamma_k-1$.  The proof is essentially the same as in the case that $|{\mathcal I} | = 1$, so we give the proof in the degree 1 setting to avoid computation.  Hence we assume that $aX^m \in I$ for $m \geq 1$ and show that $aX^{m-1} \in I$.  (We  are omitting the subscripts on $X$ and $\zeta$,  because of the
reduction to the $|{\mathcal I} | = 1$ case.)

By the definition of $I$, it is clear that $a(X-\zeta)^m \in I$.  Then it follows that $aX^m - a(X-\zeta)^m = a(X^m-(X-\zeta)^m) \in I$, and after simplification using the identity 
$$x^m-y^m = (x-y)(x^{m-1}+x^{m-2}y+  \cdots +xy^{m-2} + y^{m-1}),$$ 
we have that 
$$\zeta a(X^{m-1} + X^{m-2}(X-\zeta) + X^{m-3}(X-\zeta)^2 + \cdots + (X-\zeta)^{m-1}) \in I.$$
Since $\zeta a X^{m-i}(X-\zeta)^{i-1} \in I$ for all $i \ge 2$, it follows that $\zeta a X^{m-1} \in I$, and thus $a X^{m-1} \in I$.  By induction on $m$, we may conclude that $a \in I$.
\end{proof}

\begin{lem}\label{lem:divisAlgAnalogue}
If $\delta \in \Gamma$, then $X^\delta \in R + \sum_{i \in \mathcal I} A(X_i - \zeta_i)$.  
\end{lem}
\begin{proof}  The proof is by induction on $| \delta|$.  If $\delta_i = 0$ for all $i$, then $X^\delta = 1 \in R$.   So we suppose that $\delta_k>0$.  Then 
\begin{align*}
X^\delta &= \Big(\prod_{i \in {\mathcal I}, i \neq k} X_i^{\delta_i}\Big)X_k^{\delta_k-1} (X_k - \zeta_k  + \zeta_k) \\
&= X^{\delta'} (X_k - \zeta_k) + \zeta_k X^{\delta'},
\end{align*}
where $\delta_i' = \delta_i$ for $i \neq k$, and $\delta_k' = \delta_k-1$.   Since $X^{\delta'} (X_k- \zeta_k)  \in \sum_{i \in \mathcal I} A(X_i - \zeta_i)$,  and $|\delta'| < | \delta|$, we have by induction that $X^\delta \in R + \sum_{i \in \mathcal I} A(X_i - \zeta_i)$.
\end{proof}

\begin{rem}
It is evident from the proof of Lemma \ref{lem:divisAlgAnalogue} that in fact $X^\delta \in \kk 1 + \sum_{i\in \mathcal I} A(X_i - \zeta_i)$ for all $\delta \in \Gamma$.  In applying the lemma, however, we only need that $X^\delta \in R + \sum_{i\in \mathcal I} A(X_i - \zeta_i)$.
\end{rem}

\begin{thm} \label{thm:ann_A(w)}
Suppose that $V = Aw$ is a Whittaker module of type $\zeta$ for $A = R(\phi,t)$, and let $Q = \ann_R(w)$.  Then $\ann_A(w) = AQ + \sum_{i \in\mathcal I} A(X_i-\zeta_i)$.
\end{thm}
\begin{proof}
Let $a \in \ann_A(w)$.    As in the proof of Lemma \ref{lem:Rfree}  (in particular, as in \eqref{eq:Xpoly}),  there exists  $\gamma \in \Gamma$ such that $aX^\gamma = \sum_{\delta \in \Gamma} r_\delta X^\delta$ for $r_\delta \in R$.  Since $aX^\gamma w = \prod_i \zeta_i^{\gamma_i} aw = 0$, it follows that $aX^\gamma \in \ann_A(w)$.   Lemma \ref{lem:divisAlgAnalogue} implies  that $aX^\gamma = \sum_{\delta \in \Gamma} r_\delta X^\delta \in R + \sum_{i \in \mathcal I} A(X_i - \zeta_i)$.  Thus we may write $aX^\gamma = r + b$, with $r \in R$ and $b \in \sum_{i \in \mathcal I} A(X_i-\zeta_i) \in \ann_A(w)$.  Since $aX^\gamma \in \ann_A(w)$ and $b \in \ann_A(w)$, it must be that $r \in \ann_A(w) \cap R = Q$.  Thus $aX^\gamma \in AQ + \sum_{i=1}^n A(X_i-\zeta_i)$, and by Lemma \ref{lem:aXm}, we have $a \in AQ + \sum_{i \in \mathcal I}  A(X_i-\zeta_i)$.  The other containment is clear.
\end{proof}

\begin{cor}
If $V = Aw$ is a Whittaker module of type $\zeta$ for $A = R(\phi,t)$, then 
$$V \cong A \Bigg / \left( AQ + \sum_{i \in \mathcal I} A(X_i-\zeta_i) \right),$$
where $Q = \ann_R(w)$.
\end{cor}

Since $\ann_A(V) \subseteq \ann_A(w)$, we have the next corollary.

\begin{cor}
If $V = Aw$ is a Whittaker module of type $\zeta$ for $A = R(\phi,t)$, then $\ann_A(V) \subseteq AQ + \sum_{i\in \mathcal I} A(X_i-\zeta_i)$.
\end{cor}

If $V=Aw$ is a one-dimensional Whittaker module, then $\ann_A(V) = \ann_A(w)$, which
implies the following result.

\begin{cor}
Suppose that $V = Aw = \kk w$ is a one-dimensional Whittaker module of type $\zeta$ for $A = R(\phi,t)$.  Then $\ann_A(V) = AQ + \sum_{i\in \mathcal I} A(X_i-\zeta_i)$, and there exists
an  $\kk$-algebra homomorphism  $\theta: R \to \kk$ such that  $r w = \theta(r)w$
for all $r \in R$ and $\ker \theta = Q$.  
\end{cor}

\section{Whittaker vectors}

Assume that $(V,w)$ is a Whittaker pair of type $\zeta$ for the generalized
Weyl algebra  $A = R(\phi,t)$.
Let $\hbox{\rm Wh}_\eta(V)$ denote the set of all Whittaker vectors of type $\eta = (\eta_i)_{i \in \mathcal I}$
in $V$.   In this section we describe how Whittaker vectors are related
to eigenvalues of the automorphisms $\phi_i$ and how they  can be used
to deduce information about the module $V$.  We note that Lemma \ref{lem:Whvect} and Corollary \ref{cor:Whvectu} are true even if $R$ is noncommutative.  For the remaining results in this section, we must assume that $R$ is commutative.  The next result is apparent.

\begin{lem}\label{lem:Whvect}  Let $(V,w)$ be a Whittaker pair of type $\zeta$ with
$Q = \ann_R(w)$.  Then for  $v = rw \in V$, the following are equivalent:   \begin{itemize}
\item[\rm (a)] $X_iv = \eta_i v$, for all $i \in \mathcal I$;
\item[\rm (b)] $\eta_i r w = X_i r w = \phi_i(r)\zeta_i w$, for all $i \in \mathcal I$;
\item[\rm (c)] $\phi_i(r) \equiv \zeta_i^{-1}\eta_i r  \mod Q$;
\item[\rm (d)]  $r+Q$ is an eigenvector for the induced linear transformation $\overline {\phi_i}$ on $R/Q$ with eigenvalue $\zeta_i^{-1}\eta_i$ for all $i \in \mathcal I$.
\end{itemize} 
If $\eta_i \neq 0$ for all $i \in {\mathcal I}$, then $v \in \hbox{\rm Wh}_\eta(V)$.  \end{lem} 

\begin{cor} \label{cor:Whvectu}   For the universal Whittaker pair $(V_{\mathfrak u}, w_{\mathfrak u})$ of
type $\zeta$,
$0 \neq rw_{\mathfrak u} \in \hbox{\rm Wh}_\eta(V_{\mathfrak u})$ if and only
if $r$ is an eigenvector of $\phi_i$ with eigenvalue $\zeta_i^{-1}\eta_i$ for all $i \in \mathcal I$.  
\end{cor}

\begin{proof}  This follows directly from  Lemma \ref{lem:Whvect} and
the fact that $\ann_R(w_{\mathfrak u}) = 0$ (see Lemma \ref{lem:universalWhitt}).  
\end{proof}  

\begin{prop}\label{prop:rad}  Assume $R$ is commutative, and let  $(V,w)$ be a Whittaker pair  with 
$w \neq 0$.   Set $Q = \ann_R(w)$ and let
$$P = \sqrt{Q} = \{ r \in R \mid r^k \in Q \ \hbox{\rm for some} \ k \}.$$
Then $Pw$ is a submodule of $V$ and $Pw \neq V$.   
\end{prop} 

\begin{proof}   It suffices to note that  $P$ is a $\phi$-stable ideal of $R$ containing $Q$, and $P \neq R$  since $1 \not \in P$, and the rest follows from Lemmas \ref{lem:submodulesOfR} and \ref{lem:Whittquo}.  \end{proof}

\begin{cor}  Assume $R$ is commutative, and let  $(V,w)$ be a Whittaker pair of type $\zeta$ with $Q = \ann_R(w)$.
Assume $0 \neq v \in \hbox{\rm Wh}_\eta(V)$ and $\lambda_i:= \zeta_i ^{-1}\eta_i $ is not
a root of unity for some $i \in \mathcal I$.   If $V$ is simple, then $V$ is infinite-dimensional.  \end{cor}

\begin{proof}  Since $V$ is simple, we know by Proposition \ref{prop:rad}  that $ \sqrt{Q}w = 0$,  and hence
that $ \sqrt{Q} = Q$.     Thus,  if $r w \neq 0$, then $r^k w \neq 0$ for all $k \geq 1$.  
Now suppose that $v = rw$ is a nonzero Whittaker vector of type $\eta$ 
and set $\lambda_i =  \zeta_i ^{-1}\eta_i $.     By induction and Lemma \ref{lem:Whvect}   it follows  that 
   \begin{equation}\label{eq:rk} X_i r^k w = \lambda_i^k \zeta_i r^k w\end{equation} for all $k \geq 1$.     
Relation \eqref{eq:rk}  implies  that $r^k w$ is a nonzero eigenvector for $X_i$ with eigenvalue $\lambda_i^k \zeta_i$.       As these values are all distinct because $\lambda_i$ is not a root of unity, the vectors $r^k w$ for $k \geq 1$ must be linearly independent.  Thus, $V$ is infinite-dimensional.    \end{proof}

\begin{rem}
It is evident from the proof of the previous result that if we replace the assumption that $V$ is simple with the assumption that $Q = \sqrt Q$, the conclusion remains true.
\end{rem}

For the remainder of the section  we assume that $(V,w)$ is a fixed Whittaker pair of type $\zeta$ with $Q = \ann_R(w)$  for the generalized Weyl algebra  $A = R(\phi,t)$, where $R$ is commutative,  
and we set  
 \begin{eqnarray}\label{eqn:Whs}  S  &=& \setof{s \in R}{X_isw = \zeta_i sw \  \hbox{\rm for all}\ i \in {\mathcal I}} \\  &=& 
  \setof{s \in R}{s - \phi_i(s) \in Q \ \hbox{\rm for all}\ i \in \mathcal I}. \nonumber \end{eqnarray}
Note the second equality comes from Lemma \ref{lem:Whvect}, and
$Sw  =  {\rm Wh}_\zeta(V)$.

\begin{lem} If $R$ is commutative, and $S$ is as in \eqref{eqn:Whs}, then $S$ is a subring of $R$ and 
$Q = \ann_R(w)$  is an ideal of $S$.   
\end{lem}
\begin{proof}
If $s_1, s_2 \in S,$ then $X_is_1s_2w = \phi_i(s_1) X_is_2w = \phi (s_1) \zeta_i s_2w = \zeta_is_2 \phi (s_1) w$ \ $ = \zeta_i s_2s_1w = \zeta_i s_1s_2w.$
\end{proof}

\begin{lem}
\label{lem:Aendos} Assume $R$ is commutative and $(V,w)$ is a Whittaker pair of type $\zeta$, and let 
 $\pi : A \to {\rm End} (V)$ be the corresponding representation of $A$.  Then for
 $S$ as in \eqref{eqn:Whs}, \  $\pi (S)= {\rm End}_A (V)$.   
\end{lem}
\begin{proof}
It is clear that $srv = rsv$ for  $s \in S, r \in R$,  and $v \in V.$  We must show that $s X_iv = X_isv$ and $sY_iv=Y_isv$ whenever $s \in S$ and $v \in V.$   
But $X_isv = \phi_i(s) X_iv =  s X_i v$,  as
$s - \phi_i(s) \in Q \subseteq \ann_A(V).$   Similarly,  
$Y_isv = \phi_i^{-1} (s) Y_iv = sY_iv$,
since $s - \phi_i^{-1}(s)  = \phi_i^{-1}\big(\phi_i(s)-s\big)   \in Q \subseteq \ann_A (V)$.
Thus, $\pi(S) \subseteq {\rm End}_A (V)$.  

For the other direction,  let $\psi \in {\rm End}_A (V)$,  and note that $X_i \psi w = \zeta_i \psi w$,  so $\psi w \in {\rm Wh}_\zeta (V) = Sw.$  Write $\psi w = sw$ for $s \in S.$  It is easy to see that the action of $\psi$ on $V$ is determined by its action on $w$, and thus $\psi = \pi (s).$
\end{proof}

The map $S \to {\rm End}_A(V)$ defined  by $s \mapsto \pi (s)$ gives the following. 

\begin{cor}
${\rm End}_A (V) \cong S / Q$. 
\end{cor}

The following  rendition of Schur's lemma enables us to say more in the simple case.

\begin{lem} 
Suppose $\kk$ is an uncountable algebraically closed field and $A$ is a $\kk$-algebra.  If $V$ is a simple $A$-module of countable dimension over $\kk,$ then ${\rm End}_A (V) = \kk \hbox{\rm id}_V.$
\end{lem}

\begin{cor} \label{cor:whittVecVsimple}  Assume $R$ is commutative and is of countable dimension over an uncountable algebraically closed field $\kk$,  and let $(V,w)$ be a Whittaker pair of type $\zeta$ for $A = R(\phi,t)$.  
If  $V$ is simple, then  ${\rm Wh}_\zeta (V) = Sw =\mathcal Zw = \kk w$, where $\mathcal Z$ is
the center of $A$.  
\end{cor}
\begin{proof}
If $V$ is simple, ${\rm Wh}_\zeta (V) = S w = \kk w$ since $\pi(S) = {\rm End}_A(V) = \kk \hbox{\rm id}_V$ 
by Schur's Lemma.    But then $\kk w \subseteq \mathcal Zw \subseteq Sw = \kk w,$ forcing equality.
\end{proof}

\section{An expression for $\ann_A(V)$}\label{sec:Aquotient}
 
\begin{prop}\label{prop:idealIntersectR}
Assume that 
\begin{equation}\label{eqn:assumptionPhi}
\mbox{if \, $\lambda \in \Lambda \subseteq \mathbb Z^{\mathcal I}$ \, and  \, $\phi^\lambda:= \prod_{i \in \mathcal I} \phi_i^{\lambda_i} = {\rm id}_R$, \, then \, $\lambda_i = 0$ for all $i \in \mathcal I$.}
\end{equation}
(or equivalently by Proposition \ref{prop:Acent}  that the center of the generalized Weyl algebra
$A = R(\phi,t)$ is contained in $R$).  
If $R$ is commutative and $B$ is a nonzero ideal of $A$, then $B \cap R \neq 0$.    \end{prop}

\proof  Let $0 \neq  a \in  B$.   As in \eqref{eq:Xpoly}, there exists some $X^\gamma$ for $\gamma
\in \Gamma$  so that 
$$a X^{\gamma} =  \sum_{\delta \in \Gamma}  r_\delta X^\delta \in B.$$
Thus $B$ contains some nonzero polynomial in  the $X_i$ with coefficients in $R$, and we may assume   $b = \sum_{\varrho  \in \Gamma} b_\varrho X^\varrho  \in B$  is such a polynomial having the least number of nonzero terms.    Then for $r \in R$, 

$$b r = \sum_{\varrho \in \Gamma} b_\varrho  \phi^\varrho(r)X^\varrho \in B.$$

\noindent Suppose $b_{\sigma} \neq 0$.  Then the element
$$br - \phi^\sigma(r) b = \sum_{\varrho \in \Gamma} b_\varrho \big(\phi^\varrho (r) - \phi^{\sigma}(r) \big)X^\varrho \in B.$$
would have fewer nonzero terms unless $\phi^{\varrho}(r)=
\phi^\sigma(r)$ for all $\varrho$ with $b_\varrho \neq 0$.  However,  (\ref{eqn:assumptionPhi})  implies
that there  must exist an $r$ such that $\phi^\varrho(r) \neq \phi^\sigma(r)$ for  $\varrho \neq \sigma$.    Thus,  a nonzero polynomial in the $X_i$ belonging to $B$ and having a minimal number of terms has the form $sX^\sigma$ for some $s \in R$.   But then $sX^\sigma Y^\sigma$ is a nonzero element of  $B \cap R$.  \qed

\begin{cor} \label{cor:annRwNonzero}
Let $V =Rw$ be a Whittaker module for a generalized Weyl algebra $A = R(\phi, t)$ with $R$ commutative such that  (\ref{eqn:assumptionPhi}) holds.  Then $\ann_A(V) \cap R = \ann_R(w)$, so that if $\ann_A(V) \neq 0$, then $\ann_R(w) \neq 0$.
\end{cor}

Let  $A = R(\phi,t)$ be a generalized Weyl algebra with $R$ commutative,  and assume $J$ is a $\phi$-stable ideal of $R$ such that $t_i \not\in J$ for all $i \in \mathcal I$.  As before, let  $\overline r$ mean the coset $\overline r = r+J$ for all $r \in R$.  Thus, $\overline {t_i} \neq 0$.  Since $J$ is $\phi$-stable, we have the induced automorphisms $\overline {\phi_i}$ on the quotient $R/J$.  By Proposition \ref{prop:quotientGWA}, $A/AJ$ is isomorphic to the GWA $\overline A := (R/J)(\overline \phi, \overline t)$ whenever $J$ is a prime ideal of $R$.     Now if $V = Rw$ is a Whittaker module of type $\zeta = (\zeta_i)_{i \in \mathcal I}$ for $A$, and if $Q = \ann_R(w)$, then $AQ$ is always a 2-sided ideal contained in $\ann_A(V)$.  So if $Q$ is a prime ideal of $R$, we may always pass to the GWA $\overline A = (R/Q)(\overline \phi, \overline t)$ (provided $\overline {t_i}\neq 0$ for all $i$) and regard $V$ as a module for the (possibly) different GWA $\overline A$.  The annihilator $\ann_{\overline A}(V)$ may be nontrivial (if $AQ \neq \ann_A(V)$), but here is a situation where that does not happen.

\begin{thm}\label{thm:annVwhenQprime}
Let  $A = R(\phi,t)$ be a generalized Weyl algebra with $R$ commutative.  Assume $V = Rw$ is a Whittaker module of type $\zeta$, and suppose that $Q = \ann_R(w) = \ann_R(V)$ is a prime ideal such that $t_i \not\in Q$ and the induced automorphisms $\overline {\phi_i}$ on $R/Q$ satisfy  (\ref{eqn:assumptionPhi}).  Then $\ann_A(V) = AQ$.
\end{thm}
\begin{proof}
By the above considerations, we may suppose that $V$ is a Whittaker module of type $\zeta$ for the GWA $\overline A = A/AQ = (R/Q)(\overline \phi, \overline t)$.  Note that $\ann_{R/Q}(w) = \ann_{R/Q}(V) = 0$.

Consider the ideal $B = \ann_A(V) + AQ$ in $A/AQ = (R/Q)(\overline \phi, \overline t)$.  Then  $R/Q$ is a commutative domain since $Q$ is a prime ideal of $R$.  Now we have seen from Proposition \ref{prop:idealIntersectR} that when a GWA has a commutative Noetherian domain as its coefficient ring and when the automorphisms  satisfy (\ref{eqn:assumptionPhi}), then any ideal intersects the coefficient ring nontrivially.  Since $B \cap (R/Q) = 0$, it must be that $B = 0$ in $A/AQ$.  That is, $\ann_A(V) = AQ$.
\end{proof}

\section{Examples -  a brief introduction}

We apply results of the previous sections to determine the Whittaker modules for the  quantum plane and the  (quantum) Weyl algebra, and for certain generalizations of the universal enveloping algebra $U(\mathfrak{sl}_2)$  introduced by S.P. Smith and their quantum analogues.  The algebras considered here have a realization as generalized Weyl algebras $A = R(\phi,t)$,  where $R$ is commutative and $\phi$ is a single automorphism.  Because the Whittaker modules of type $\zeta$ are in bijection with the $\phi$-stable ideals of $R$, we begin by describing those ideals.   For the Smith algebras and quantum Smith algebras,  the $\phi$-stable ideals $J$ of $R$ are generated by their intersection $J \cap \mathcal Z$ with the center $\mathcal Z = \mathcal Z(A)$ of $A$.    In determining that intersection, the   description of the center of a generalized Weyl algebra in Proposition \ref{prop:Acent}  is essential.

\section{The case $\boldsymbol{R = \kk[t]}$ }\label{sec:RFt}   

An automorphism $\phi$ of the polynomial algebra $R = \kk[t]$ is necessarily given by $\phi(t) = \alpha t + \beta$ for some $\alpha,\beta \in \kk$ with $\alpha \neq 0$.  Let 
\begin{equation*}
\tilde t = (\alpha -1)t + \beta,
\end{equation*}
and note that $R = \kk [\tilde t]$ as long as $\alpha \neq 1$.  Since $\phi (\tilde t) = \alpha \tilde t$, it is evident that $\phi^\ell = \hbox{\rm id}_R$ if $\alpha \neq 1$ is a primitive $\ell$th root of unity, and $\phi$ has infinite order if $\alpha$ is not a root of unity.

\begin{lem}\label{lem:phiStabR=Kt}
Assume $\alpha \neq 1$ and $J$ is a nonzero proper $\phi$-stable ideal of $R = \kk [t] = \kk [\tilde t]$, and let $f(\tilde t)$ be the unique monic generator of $J$.  If $\alpha$ is not a root of unity, then there exists $n > 0$ such that $f(\tilde t) = \tilde t^n$.  If $\alpha$ is a primitive $\ell$th root of unity, then there exist $n \ge 0$ and scalars $c_k \in \kk$ such that $f(\tilde t) = \tilde t^n \sum_{k \ge 0} c_k\tilde t^{k\ell}$.   
\end{lem}
\begin{proof}
Since $\phi (\tilde t) = \alpha \tilde t$, it is clear that polynomials of the stated forms (in either case) generate $\phi$-stable ideals.

Conversely, if $J$ is a $\phi$-stable ideal of $\kk [\tilde t]$, let  $f = f(\tilde t) = \sum_{j \ge 0}^n  a_j \tilde t^j$, $a_n = 1$,  be the unique monic polynomial in $\tilde t$ (of minimal degree) generating $J$.  
Then $\alpha^n f- \phi( f) = \sum_{j \ge 0}^n  a_j \left(\alpha^n-\alpha^j\right) \tilde t^j \in J$.  
  If $\alpha$ is not a root of unity, then $a_j = 0$ for all $j \neq n$ and $f(\tilde t) = \tilde t^n$.  Suppose $\alpha$ is a primitive $\ell$th root of unity.   Since $\alpha^{n-j} =1$ must hold whenever  $a_j \neq 0$,  we have  $j \equiv n$ ${\rm mod} \, \ell$ for each such $j$, so the polynomial $f$ has the desired form. 
\end{proof}

\begin{cor}\label{cor:maxStabR=KtRoot}  Suppose $\kk$ is algebraically closed, and  $J$ is a  $\phi$-stable ideal of $R = \kk [\tilde t]$, where $\phi(\tilde t) = \alpha \tilde t$ and  $\alpha \neq 1$.   Let $f(\tilde t)$ be the unique monic generator of $J$.  Then $J$ is a proper maximal $\phi$-stable ideal if and only if  $f(\tilde t) = \tilde t$ when $\alpha$ is not a root of unity, and  $f(\tilde t) = \tilde t$ or $f(\tilde t) = \tilde t^\ell - \xi$ for some nonzero $ \xi \in \kk$  when $\alpha$ is a primitive $\ell$th root of unity.  
\end{cor}

Next we examine in detail the Whittaker modules for the generalized Weyl algebra $A = R(\phi,t)$ constructed from $R = \kk[t]$ and the automorphism $\phi$.  Thus $YX = t$, $XY = \phi(t) = \alpha t + \beta$, and  as before, we assume $\alpha \neq 1$.    Let $(V,w)$ be a Whittaker pair of type $\zeta$ for $A$, and let $Q = \ann_R(w)$.   By Theorem \ref{thm:bijectIdealsModules}, if $Q = 0$, then $V$ is isomorphic to the universal Whittaker module $V_{\mathfrak u}$ of  type $\zeta$, so we will assume $Q \neq 0$.

\begin{subsec}  {$\boldsymbol{\alpha}$ \textbf{is not a root of unity}}    \end{subsec}

When  $\alpha$ is not a root of unity, then $Q = R \tilde t^n$ for some $n \geq 1$.  Since $V = Rw \cong  V_Q = R/Q$, it is clear that 
$\setof{v_k: = \tilde t^k w}{0 \le k\le n-1}$ is a basis  of $V$ and $\dim_\kk V = n$. The action of $A$ on $V$ is given as follows
\begin{eqnarray} \label{eq:WR=Kt} 
&& t v_k=  (\alpha-1)^{-1}
\big(v_{k+1} - \beta v_k\big),   \\
&& X v_k= \alpha^k \zeta v_k,     \qquad Y v_k =  \zeta^{-1}  \alpha^{-k} (\alpha-1)^{-1}
\big(v_{k+1} - \beta v_k\big), 
 \nonumber \end{eqnarray}

\noindent where $v_{n} = 0$.   Since submodules of $V$ correspond to $\phi$-stable ideals of $R$ containing $Q$, the submodules of $V$ are given by $V = V_0 \supseteq V_1 \supseteq \cdots \supseteq V_n = 0$, where $V_k = Av_k = Rv_k$  is a Whittaker (sub)module with cyclic Whittaker vector $v_k$ of type $\alpha^k \zeta$, and  $\{v_k, v_{k+1}, \ldots, v_{n-1} \}$ is a basis for $V_k$.

\begin{thm} \label{thm:ann-1varPoly}
Let $R = \kk [t]$, and let $\phi : R \to R$ be the algebra automorphism given by  
$\phi (t) = \alpha t + \beta$, where $\alpha$ is not a root of unity.  Let $(V,w)$ be a Whittaker pair of type $\zeta$ for $A = R(\phi,t)$, and assume that $Q := \ann_R(w) =R\tilde t^n$ for some $n \geq 1$, where $\tilde t = (\alpha-1)t + \beta$.  Then  $V$ has a basis $v_k = \tilde t^kw$, $k=0,1,\dots, n-1$, and the action of $A$ on $V$ is given by \eqref{eq:WR=Kt}. Moreover, 
$$\ann_A(V) = \sum_{j=0}^n A {\tilde t}^{n-j}\Bigg( \prod_{k=0}^{j-1} \left(X - \alpha^k \zeta \right)\Bigg).$$
If $V$ is simple, then $n = 1$,  $\dim_\kk V = 1$, and $\ann_A(V) = A \tilde t + A (X -\zeta)$.  Any $n$-dimensional space $V$ with a basis $v_k$, $k=0,1,\dots, n-1$, and $A$-action given by  \eqref{eq:WR=Kt} is a Whittaker module of type $\zeta$ with cyclic Whittaker vector $w = v_0$ and with $Q = R \tilde t^n$.
\end{thm}
\begin{proof} All that remains to be shown is that $\ann_A(V)$ equals the expression on the right.   It is straightforward to verify that any element of the stated form annihilates every basis vector $\tilde t^kw$, $k = 0,1 \ldots, n-1$, and therefore annihilates $V$.

For the other inclusion, we proceed by induction on $\dim V$.  Let 
$$K = \sum_{j=0}^n A {\tilde t}^{n-j}\Bigg( \prod_{k=0}^{j-1} \left(X - \alpha^k \zeta \right)\Bigg),$$
and suppose that $a \in \ann_A(V)$.  Since $\ann_A(V) \subseteq \ann_A(w) = AQ + A (X - \zeta)$, we may write $a = a_1 q + a_2(X - \zeta)$ with $q \in Q \subseteq K$ and $a_1,a_2 \in A$.  Both $a$ and $a_1q$ annihilate $V$, so $a_2(X - \zeta)$ must annihilate $V$.  In particular, $a_2(X - \zeta)$ must annihilate the submodule $V_1 = A \tilde t^1 w \subseteq V$.  Since $V_1 = {\rm span}_\kk \{ \tilde t^k w, \mid k=1,\dots,n-1\}$ and $(X - \zeta) \tilde t^k w = \zeta (\alpha^k -1) \tilde t^kw$, it follows that $a_2$ annihilates the Whittaker module $V_1$ (of type $\alpha \zeta$).  Since $\ann_R(\tilde t^1 w) =   R\tilde t^{n-1}$, we may claim by induction that 
\begin{equation*}
a_2 \in \sum_{j=0}^{n-1} A \tilde t^{n-1-j} \Bigg( \prod_{k=0}^{j-1} \left(X - \alpha^k\alpha \zeta \right)\Bigg)  = \sum_{j=1}^n A \tilde t^{n-j} \Bigg( \prod_{k=1}^{j} \left(X - \alpha^k \zeta \right)\Bigg).
\end{equation*}
Hence  \ $a = a_1 q + a_2(X - \zeta) \in K$, 
as desired.
\end{proof}
\medskip

\begin{subsec}  {$\boldsymbol{\alpha}$ \textbf{is a root of unity, \  $\boldsymbol{\alpha} \neq 1$}}    \end{subsec}

Now suppose that $\alpha$ is a primitive $\ell$th root of unity.   We have observed earlier that $\phi^\ell = \hbox{\rm id}_R$ in this case and that the center $\mathcal Z$  of $A$ is generated by $X^\ell$, $Y^\ell$ and the set $R^\phi$ of elements of $R$ fixed by $\phi$, which are the polynomials in $\kk[\tilde t^\ell]$.  We will assume that $\kk$ is algebraically closed, and the Whittaker module $V$ is simple.  
Then since $Q = \ann_R(w)$ is a maximal $\phi$-stable ideal, it follows from Corollary \ref{cor:maxStabR=KtRoot}, that $Q  = R( \tilde t^\ell - \vartheta^\ell)$ for some nonzero $\vartheta \in \kk$ or
$Q = R \tilde t$.    In the former case,    $V \cong V_{\mathfrak u}/QV_{\mathfrak u} = R/Q$ is $\ell$-dimensional, and the vectors $u_k: =  \vartheta^{-k} \tilde t^k w,  k=0,1,\dots, \ell-1$, 
determine a basis for $V$.   Moreover, from   (\ref{eqn:univModueXaction}) and 
 (\ref{eqn:univModueYaction}) we see   that
$X^\ell = \zeta^\ell \hbox{\rm id}_V$,  $Y^\ell =  
\zeta^{-\ell} (\alpha-1)^{-\ell} \alpha^{-(\ell -1) \ell /2} (\vartheta^\ell - \beta^\ell)\hbox{\rm id}_V$
 and the following hold:
 
\begin{eqnarray}\label{eq:WRKtunit}  
&& X u_k = \zeta \alpha^k u_k,   \quad  Y u_k = \zeta^{-1}\alpha^{-k} (\alpha-1)^{-1}(\vartheta u_{k+1}-\beta u_k),   \\
&& t u_k = (\alpha-1)^{-1}(\vartheta u_{k+1}-\beta u_k), \ \    \nonumber
\end{eqnarray}
where subscripts should be read mod $\ell$.  Now when $Q = R\tilde t$, then $V \cong R/Q = \kk 1$, 
and  $V = \kk w$.    In this case,  $tw = -(\alpha-1)^{-1} \beta w$,  $Xw = \zeta w$, and $Y w = -\zeta^{-1}(\alpha-1)^{-1}\beta w$.    In summary we have  
\begin{thm} \label{thm:ann-1varPolyunity}
Let $R = \kk [t]$  where $\kk$ is an algebraically closed field, and let $\phi : R \to R$ be the algebra automorphism given by  $\phi (t) = \alpha t + \beta$, where $\alpha$ is a primitive $\ell$th  root of unity.  Assume $V$ is a simple Whittaker module of type $\zeta$ for $A$.  Then either
\begin{itemize}
\item[\rm {(i)}] $\dim V = \ell$,   $Q = \ann_R(w) = R( \tilde t^\ell - \vartheta^\ell)$ for some $\vartheta \neq 0$, and $V$ has a basis $u_k, k = 0,1, \dots, \ell-1$, so that the action of $A$ on $V$ is given by \eqref{eq:WRKtunit};   or  
\item[\rm {(ii)}] $\dim V = 1$,  $Q = \ann_R(w) = R \tilde t$, and $V = \kk w$, where
\begin{equation} \label{eq:onedimru} tw = -(\alpha-1)^{-1} \beta w, \quad Xw = \zeta w, \quad  
Y w = -\zeta^{-1}(\alpha-1)^{-1}\beta w. \end{equation}  \end{itemize}     
\noindent Conversely, any $\kk$-vector space $V$ of dimension $\ell$ (or $1$)  having an $A$-action given by 
\eqref{eq:WRKtunit} (or  by \eqref{eq:onedimru})  determines a simple Whittaker module of type $\zeta$   
for $\zeta \neq 0$.  \end{thm}

\begin{subsec}  {$\boldsymbol{\alpha = 1}$}    \end{subsec}

If $\alpha  = 1$ and $\beta = 0$, the $\phi$-stable ideals are just ordinary ideals $J$ and every $V = R/J$ with  $A$-action inherited from \eqref{eqn:univModueRaction}-\eqref{eqn:univModueYaction} is a Whittaker module of type $\zeta$.     If $\alpha = 1$ and $\beta \neq 0$, it is evident that $R$ contains no nontrivial proper $\phi$-stable ideals when $\kk$ has characteristic 0.  Thus there is, up to isomorphism, only one Whittaker module, namely the universal one $V_{\mathfrak u} = R$, and it is necessarily simple.  In particular, when  $\beta = -1$, the algebra $A = R(\phi, t)$ is the Weyl algebra $A_1$, and \eqref{eq:Weylex}  gives the $A$-action on $V_{\mathfrak u}$ 
in this special case.
 
When the characteristic of $\kk$ is $p>2$ and $\phi(t) = t+\beta$, where $\beta \neq 0$,  then $\phi^p(t) = t + p\beta  = t$, and $\phi$ has order $p$.   In this case, it follows from Proposition \ref{prop:Acent}   that the center $\mathcal Z$  of $A = R(\phi,t)$ is generated by $X^p, Y^p$, and the set $R^\phi$ of elements of $R$ fixed by $\phi$.  It is straightforward to verify that $t^p - \beta^{p-1} t$ is fixed by $\phi$.  Define 
$$ z_n = \begin{cases}   t^n  & \quad \hbox{\rm if} \ n \not \equiv 0  \mod p \\
(t^p - \beta^{p-1} t)^{n/p}  & \quad \hbox{\rm if} \ n  \equiv 0  \mod p \end{cases}$$
so that the set $\setof{z_n}{n \in \z_{\geq 0}}$ is a basis for $R$.    Now let $g = \sum_{k=0}^m  g_k z_k$, where $g_0, \ldots, g_m \in \kk$, and assume that $g$ is fixed by $\phi$.  It can be shown that $g - \phi(g) \neq 0$ unless $g = \sum_{k \equiv 0 \mod p}  g_k z_k$, which is a polynomial in $t^p - \beta^{p-1} t$.   Thus,  when $A$ has characteristic $p > 2$,   the center   of $A$ is generated by $X^p, Y^p, t^p-\beta^{p-1}t$.  Observe that $z_p^j = z_{jp}$, so that $R^\phi = R \cap \mathcal Z$ is the polynomial algebra $\kk[z_p]$, and $R$ is a free $R^\phi$-module with basis $1,t, \dots, t^{p-1}$.

We claim that  if $J$ is a $\phi$-stable ideal of $R$, then  $J = R(J \cap \mathcal Z)$.   If this assertion is false,  then  there  is a polynomial $f =  \sum_{j=0}^n s_j t^j \in J \setminus R(J \cap \mathcal Z)$ with coefficients in $R^\phi$ of least degree in $t$.  Thus, $0 < n \leq p-1$, $s_n \neq 0$, and $s_j \in R^\phi$ for all $0 \leq j \leq n$.  Then $\left(\hbox{\rm id}_R - \phi\right)^n (f)  = n ! s_n \in J$.  But this implies $s_n \in J \cap \mathcal Z$ and hence that $s_n t^n \in R(J \cap \mathcal Z)$.  By minimality of $n$,  we have $f-s_nt^n \in R(J \cap \mathcal Z)$, and so $f = s_nt^n + (f-s_nt^n) \in R(J \cap \mathcal Z)$, a contradiction.   Thus, $J = R(J \cap \mathcal Z)$, so that every $\phi$-stable ideal of $R$
is centrally generated.

Theorem \ref{thm:bijectIdealsModules}  gives a bijection between isomorphism classes of Whittaker modules of type $\zeta$ and $\phi$-stable ideals of $R$  given by $V \mapsto \ann_R(w)$.  Since $Q:= \ann_R(w)$ is $\phi$-stable, we know that $Q = R(Q \cap \mathcal Z)$.

Now assume $\kk$ is algebraically closed and $V$ is a simple Whittaker module of type $\zeta$ for $A$.  Then $Q$ is a proper maximal $\phi$-stable ideal.  Since $Q \cap \mathcal Z$ is a maximal  ideal of $R^\phi = R \cap \mathcal Z = \kk[z_p]$,  we can find $\lambda \in \kk$ so that 
$Q \cap \mathcal Z =  R^\phi\big( z_p -(\lambda^p-\beta^{p-1}\lambda)\big)$.  Thus, 
$Q = R\big( z_p -(\lambda^p-\beta^{p-1}\lambda)\big)  = R\big( t^p - \beta^{p-1} t -(\lambda^p - \beta^{p-1} \lambda)\big)$, and $V \cong V_{\mathfrak u}/Q V_{\mathfrak u} = R/Q$ is $p$-dimensional.   Since $(t-\lambda)^p - \beta^{p-1}(t-\lambda) =  \prod_{k = 0}^{p-1}\big(t-(\lambda-k\beta)\big)$,
the eigenvalues of $t$ on $V$ are of the form  $\lambda-k\beta$ for $k=0,1,\dots, p-1$.  
Let  $v_0 = \prod_{k = 1}^{p-1} \big(t-(\lambda-k\beta)\big)$, and observe
that $(t-\lambda)v_0 = 0$.     Since $X^p$ is central, it acts as a scalar on $V$, and 
{f}rom $X^p w = \zeta^p w$ we see that scalar is $\zeta^p$.   
The vectors  $v_k := \zeta^{-k} X^k v_0$ for $k=0,1,\dots, p-1$
determine a basis for $V$  and relative to this basis, the $A$-action is given by \begin{eqnarray}\label{eq:pactionalone}  t v_k &=& (\lambda-k\beta)v_k    \\
X v_k &=&  \zeta v_{k+1}  \quad (\hbox{\rm subscripts mod} \  p )   \nonumber \\
Y v_k &=& \zeta^{-1}( \lambda-(k-1)\beta)v_{k-1} 
 \quad  (\hbox{\rm subscripts mod} \  p ).  \nonumber
\end{eqnarray}
In particular,    $Y^p  = \zeta^{-p}\prod_{k=0}^{p-1}(\lambda -k\beta) = \zeta^{-p}(\lambda^p-\beta^{p-1}\lambda)\hbox{\rm id}_V$.      To summarize,
we have

\begin{thm}\label{thm:palone}   Assume   $\kk$ is algebraically closed of characteristic $p > 2$, and let $\phi : R \to R$ be the algebra automorphism of $R = \kk[t]$ given by  $\phi (t) = t + \beta$ for  $\beta \neq 0$.    If $V$ is a simple Whittaker module for $A$ of type $\zeta$, then $V$ has dimension $p$, and there is a basis $v_k, k = 0,1,\dots, p-1$,  so that the action of $A$ on $V$ is given by \eqref{eq:pactionalone} for some scalar $\lambda \in \kk$.  The vector $w = v_0+v_1 + \cdots  + v_{p-1}$ is a cyclic Whittaker vector of type $\zeta$ for $V$.   Moreover,  $Q = \ann_R(w) = R\big( t^p-\beta^{p-1}t - (\lambda^p-\beta^{p-1}\lambda)\big)$, and $\ann_A(w) = AQ + A(X-\zeta)$.    \end{thm}

\begin{rem} {We have shown that  when $\kk$ has  characteristic $p > 2$ and $\phi : R \to R$ is the automorphism of $R = \kk[t]$ given by  $\phi (t) = t + \beta$ for  $\beta \neq 0$, then the $\phi$-stable ideals of $R$ are centrally generated.  The hypotheses of Theorem \ref{thm:Whittsimple} are satisfied, and so $V_{\mathfrak u, \mathcal K} = V_{\mathfrak u}/\mathcal K V_{\mathfrak u}$  is a simple Whittaker module for every maximal ideal   $\mathcal K$  of $\mathcal Z$.    Thus, when $\kk$ is algebraically closed, $V_{\mathfrak u,\mathcal K}$ has dimension $p$ by Theorem \ref{thm:palone}, and there is a basis $v_k, k = 0,1,\dots, p-1$,  so that the action of $A$ on $V$ is given by \eqref{eq:pactionalone} for some scalar $\lambda \in \kk$.   The ideal $Q =  \ann_R(w)$ is generated by the element  $t^p-\beta^{p-1}t - (\lambda^p-\beta^{p-1}\lambda) \in Q \cap \mathcal Z = R^\phi \cap \mathcal K$.   Thus,  $V_{\mathfrak u, \mathcal K} \cong V_{\mathfrak u}/Q V_{\mathfrak u} = R/Q$,  where $Q = R\big( t^p-\beta^{p-1}t - (\lambda^p-\beta^{p-1}\lambda)\big)$. }
\end{rem}

Next we consider some well-known generalized Weyl algebras which fit into the pattern of arising from the polynomial algebra  $R = \kk[t]$.

\begin{subsec}{\textbf{The quantum plane: \quad  $\boldsymbol{R = \kk [t]}$ and $\boldsymbol{\phi (t) = \alpha  t}$ for $\boldsymbol{\alpha \neq 0, 1}$}}
\end{subsec}

In the generalized Weyl algebra $A = R(\phi, t)$ constructed from the data $R = \kk [t]$ and $\phi (t) = \alpha  t$, we have $YX = t$, $XY = \alpha t$ so that $A$ is a quantum plane.    When $\alpha$ is not a root of unity, the simple Whittaker modules are one-dimensional, $V = \kk w$,  with the action of $A$ given by
$$X w = \zeta w, \quad Y w =  0, \quad  tw = 0, $$
and  $\ann_A(V) = A\tilde t + A(X - \zeta)$, where $\tilde t = (\alpha-1)t$.
 
When $\alpha$ is a primitive $\ell$th root of unity and $\kk$ is algebraically closed, the simple Whittaker modules of type $\zeta$
are $\ell$-dimensional with  basis $u_k$, $k=0,1,\dots, \ell-1$,  and $A$-action given by
\begin{eqnarray*}\label{eq:WRKtunity}  
&& X u_k = \zeta \alpha^k u_k,   \quad  Y u_k = \zeta^{-1}\alpha^{-k}  (\alpha-1)^{-1}\vartheta u_{k+1}  \\
&& t u_k = (\alpha-1)^{-1}\vartheta u_{k+1}, \quad (\hbox{\rm subscripts mod} \ \ell) \ \ 
\end{eqnarray*}
for some scalar $\vartheta \neq 0$, or they
are one-dimensional $V = \kk w$ with
$Xw = \zeta w$, $Y w = 0$, and $t w = 0$.      
In the first case  $\ann_A(w) = A(\tilde t^\ell-\vartheta) + A(X-\zeta)$, while in the
second,  $\ann_A(V) = A\tilde t + A(X - \zeta)$.

\begin{subsec}\label{subsec:quantumWeyl}{\textbf{The quantum Weyl algebra $\boldsymbol{A_{q,1}}$:  \\  $\boldsymbol{\qquad \ R = \kk [t]}$ and $\boldsymbol{\phi (t) = q^{-1}(t-1)}$}}\end{subsec}

Fix $q \in \kk^\times$.  Let $R = \kk [t],$ and define $\phi : R \to R$ by $\phi (t) = q^{-1}(t-1).$  The algebra $A = R(\phi, t)$ is commonly referred to as the quantum Weyl algebra and is often denoted $A_{q, 1}$.  We may view $A$ as the unital algebra generated by elements $X$ and $Y$ over the field $\kk$ with relations $YX - qXY = 1$.    In the special case that $q = 1$, we obtain the (first) Weyl algebra.    In terms of the notation $\phi (t) = \alpha t + \beta$ from the previous section,   we have $\alpha = q^{-1}$ and $\beta = -q^{-1}$.
  
When $q$ is not a root of unity,  then by Theorem \ref{thm:ann-1varPoly},  the simple Whittaker modules of type $\zeta$   are one-dimensional,  $V = \kk w$, with  the action of $A$ given by
\begin{equation}\label{eq:qWeyl} X w = \zeta w, \quad Y w =  \zeta^{-1}(1-q)^{-1}w, \quad  tw = (1-q)^{-1}w,  \end{equation} 
and $\ann_A(V) = A \tilde t + A (X -\zeta)$,  where $\tilde t = (q^{-1}-1)t - q^{-1} = q^{-1}\big((1-q)t - 1\big)$.

When $q$ is a primitive $\ell$th root of unity for $\ell \geq 2$,  and $\kk$ is algebraically closed, Theorem \ref{thm:ann-1varPolyunity} implies that the simple Whittaker modules of type $\zeta$ are $\ell$-dimensional with  a basis $u_0,u_1, \dots, u_{\ell-1}$ and $A$-action given by 
\begin{eqnarray}\label{eq:qWeyll}  
&& X u_k = \zeta q^{-k} u_k,   \quad  Y u_k = \zeta^{-1}q^{k+1} (1-q)^{-1}(\vartheta u_{k+1}+q^{-1} u_k), \quad \quad  \\
&& t u_k = q(1-q)^{-1}\vartheta u_{k+1}+(1-q)^{-1} u_k, \ \    \nonumber
\end{eqnarray}
where subscripts should be read mod $\ell$ and $\vartheta \neq 0$,   or they are one-dimensional, $V = \kk w$, with the $A$-action 
\begin{equation} \label{eq:onedim} tw = (1-q)^{-1} w, \quad Xw = \zeta w, \quad  
 Y w = \zeta^{-1}(1-q)^{-1}w. \end{equation} 
 In the first case $\ann_A(w) = A (\tilde t^\ell-\vartheta) + A (X -\zeta)$, and in the second $\ann_A(V) = A \tilde t + A (X -\zeta)$,  where $\tilde t = q^{-1}\big((1-q)t - 1\big)$.

When $q = 1$,  then $A = R(\phi,t)$ is the (first) Weyl algebra $A_1$,  and
$\phi(t) = t-1$.   As we have discussed earlier, when  $\kk$ has characteristic 0, the universal Whittaker module $V_{\mathfrak u} = R$ of type $\zeta$  is simple 
and the $A$-action is given by \eqref{eq:Weylex}.  
When $q =1$ and  $\kk$ has characteristic $p > 2$, we may apply Theorem \ref{thm:palone}  with
$\beta = -1$ to deduce the following.  

\begin{thm}\label{thm:Weylp}  
Assume $\kk$ is algebraically closed of characteristic $p > 2$, and let $A$ be the (first) Weyl algebra over $\kk$ so that  $YX = t$, $XY = t-1$, and $YX-XY = 1$.  If $V$ is a simple Whittaker module for $A$ of type $\zeta$, then $V$ has dimension $p$, and there is a basis $v_k, k = 0,1,\dots, p-1$,  so that the action of $A$ on $V$ is given by 
$$t v_k = (\lambda+k)v_k,  \quad X v_k =  \zeta v_{k+1}, \quad   
Y v_k  =  \zeta^{-1}( \lambda+(k-1))v_{k-1}$$
for some scalar $\lambda \in \kk$,  (subscripts should be read  mod  $p$).  The vector $w = v_0+v_1 + \cdots  + v_{p-1}$ is a cyclic Whittaker vector of type $\zeta$ for $V$.   Moreover,  $Q = \ann_R(w) = R\big( t^p-t - (\lambda^p-\lambda)\big)$, and $\ann_A(w) = AQ + A(X-\zeta)$.    \end{thm}

\section{Smith algebras }\label{sec:Smith} 

In \cite{smith:smithAlg}, S.P. Smith introduced a family of associative algebras $\mathcal A$ which generalize the universal enveloping algebra $U(\mathfrak{sl}_2)$ of the Lie algebra $\mathfrak{sl}_2$.    These algebras are Noetherian domains with Gelfand-Kirillov dimension 3.   Smith defined a notion of weight module for the algebra $\mathcal A$ and showed there is a category of $\mathcal A$-modules analogous to the Bernstein-Gelfand-Gelfand category $\mathcal O$.   Under special assumptions, the finite-dimensional modules for $\mathcal A$  are completely reducible.  In \cite{tang:WMSmithAlg}, Tang studied Whittaker modules for the algebra $\mathcal A$ over $\mathbb C$ and obtained exact analogues of the results by Kostant in  \cite{kostant:wvrt78} for $U(\mathfrak{sl}_2)$ and by Ondrus in \cite{ondrus:Whittaker} for $U_q(\mathfrak{sl}_2)$. 

Smith's algebras $\mathcal A$ have a realization as generalized Weyl algebras, and here we show how the results we have obtained can be specialized to recover Tang's results on Whittaker modules for these algebras.   As a very special case, we obtain Kostant's results for Whittaker modules for $\mathfrak{sl}_2$.   We also apply our results to determine  the Whittaker modules
in the modular case, which was not treated in the papers of Kostant and Tang.

Fix a nonzero polynomial $s(x)$ in the algebra $\kk[x]$ of polynomials in $x$ over a field $\kk$ of characteristic not 2, and consider a unital associative algebra $\mathcal A$ over $\kk$ with generators $e,f, h$ which satisfy  the defining relations

\begin{equation}\label{eq:smith} he-eh = e,  \qquad  hf-fh = -f, \qquad  ef-fe = s(h). 
\end{equation} 
In particular, when  $s(h) = 2h$, the algebra $\mathcal A$ is isomorphic to  $U(\mathfrak{sl}_2)$.  Smith showed that there is a polynomial $r(x)$ such that

\begin{equation}\label{eq:sfromr}  s(x) = \frac{1}{2}(r(x+1)-r(x)),  \end{equation}
and the ``Casimir element,"
$$c = 2fe + r(h+1),$$
is central in $\mathcal A$.  When $\kk$ has characteristic 0, the center $\mathcal Z$  of $\mathcal A$  consists just of polynomials in $c$.  

To realize $\mathcal A$ as a GWA, let $R = \kk[h,c]$, the polynomial algebra over $\kk$ in commuting variables $h,c$,  and let $\phi$ be the automorphism of $R$ specified by $\phi(h) = h-1$, $\phi(c) = c$.    Set  $t = \frac{1}{2} (c-r(h+1))$.   Then $\phi(t) = \frac{1}{2} (c-r(h))$, and in $A = R(\phi,t)$ the following relations hold:

\begin{eqnarray*}   & YX = t = \frac{1}{2} (c-r(h+1)),  \quad  XY = \phi(t) = \frac{1}{2} (c-r(h)),     \\
& Xh = (h-1)X,  \quad  Xc = cX,   \quad   Yh = (h+1)Y,  \quad Yc = cY.
\end{eqnarray*}
Therefore, by identifying $X$ with $e$ and $Y$ with $f$,  we obtain an isomorphism between $A = R(\phi,t)$ and Smith's algebra $\mathcal A$.   In what follows we will use the GWA realization to describe the Whittaker modules for Smith's algebra.   

First suppose that $\kk$ has characteristic 0.   Then $R^\phi = R \cap \mathcal Z = \kk[c]$, and $R = \kk[h,c]$ is a free $R^\phi$-module with basis $\{h^j \mid j=0,1,\dots\}$.    Let  $J$ be a $\phi$-stable ideal of $R$.  We claim that $J = R(J \cap \mathcal Z)$.   If this is not true, there is a polynomial $f =  \sum_{j=0}^n f_j h^j \in J \setminus R(J \cap \mathcal Z)$ with coefficients in $R^\phi$ of least degree in $h$.  Thus, $n > 0$, $f_n \neq 0$, and $f_j \in R^\phi = \kk[c]$ for all $0 \leq j \leq n$.     Now  
$$\left(\hbox{\rm id}_R - \phi\right)(f)  =  \sum_{j=0}^n f_j \left(h^j - (h-1)^j \right)  = n f_n h^{n-1} + \ \hbox{\rm  lower terms in}\ h, $$
and $\left(\hbox{\rm id}_R - \phi\right)(f)   \in J$, so it follows that  $\left(\hbox{\rm id}_R - \phi\right)^n (f)  = n ! f_n  \in J$.  But this implies $f_n \in J \cap R^\phi = J \cap \mathcal Z$,  and hence that $f_n h^n \in R(J \cap \mathcal Z)$, so that $f-f_nh^n  \in J$.  The minimality of $n$ forces $f-f_nh^n \in R(J \cap \mathcal Z)$, and this gives the contradiction  $f = f_n h^n + (f-f_nh^n) \in R(J \cap \mathcal Z)$.  Thus, $J = R(J \cap \mathcal Z)$ as claimed. 

Now assume  $(V,w)$ is a Whittaker pair for $A$ of type $\zeta$ and let $Q = \ann_R(w)$.  As $Q$ is $\phi$-stable,   $Q = R (Q \cap \mathcal Z)$.     Observe that $\mathcal Z_V: = \ann_A(V)  \cap \mathcal Z  = \ann_R(V) \cap \mathcal Z = \ann_R(w) \cap \mathcal Z = Q \cap \mathcal Z$, and by Theorem \ref{thm:bijectIdealsModules},  the map  
\begin{equation}\label{annrwsmith}V \mapsto  Q = \ann_R(w)  \mapsto  Q \cap \mathcal Z = \mathcal Z_V \end{equation} 
is a bijection.   
\noindent By Theorem \ref{thm:ann_A(w)}, 
\begin{equation}\label{eq:annwsmith} \ann_A(w) = AQ + A(X-\zeta)  =  A \mathcal Z_V  + A (X-\zeta), 
\end{equation}
\noindent  which is Theorem 2.2 of \cite{tang:WMSmithAlg}.  Theorem 2.3 of \cite{tang:WMSmithAlg} establishes a one-to-one correspondence between isomorphism classes of Whitaker modules for $A$ and ideals of the center $\mathcal Z = \kk[c]$ given by $V \rightarrow \mathcal Z_V = \ann_R(w)$, as above.   (Tang assumes  $\kk = \mathbb C$,  but only characteristic 0 is necessary for these results.)

Suppose now  that $\kk$ is algebraically closed of characteristic 0, and let $V$ be a simple Whittaker module for $A$ of type $\zeta$.    Then $Q = \ann_R(w) = R(Q \cap \mathcal Z)$, where $Q \cap \mathcal Z$ is a maximal ideal
of $R^\phi = R \cap \mathcal Z = \kk[c]$.  Thus, there is $\vartheta \in \kk$ so that
$Q = R( c-\vartheta)$.   Since   $V \cong R/Q = \kk[h]$, the elements 
 $h^k w$, $k=0,1,\dots$ give a basis for $V$ and  the following hold:
\begin{eqnarray}\label{eq:actzero}  c.h^k w &=& \vartheta h^k w,  \qquad   
h.h^k w =   h^{k+1} w   \\
X.h^k w &=&\phi(h)^k \zeta w = \zeta (h-1)^k w  \nonumber  \\
Y.h^kw &=&\phi(h)^{-k} Yw = \zeta^{-1}(h+1)^k YXw = \frac{1}{2} \zeta^{-1} (h+1)^k (\vartheta - r(h+1)) w. \nonumber  
\end{eqnarray}
Now  $R/Q \cong \kk[h]$ is a domain, and the induced automorphism $\overline \phi: R/Q \to R/Q$ clearly has infinite order since $\overline \phi (h) = h-1$.  Thus if $t \not\in Q$, Theorem \ref{thm:annVwhenQprime} implies that $\ann_A(V) = A(c - \vartheta)$.  Recall that $t = \frac{1}{2} (c-r(h+1))$, where $r$ is the polynomial defined by $s(x) = \frac{1}{2}(r(x+1)-r(x))$.  If $t \in Q$, then $t - \phi (t)$ belongs to $Q$ and is a nonzero polynomial in $h$ since $t - \phi (t) = \frac{1}{2} \left( r(h) - r(h+1) \right) = \frac{1}{2} s(h) \neq 0$.   Consequently, there exists $k>0$ such that $({\rm id}_R - \phi)^k (t)$ is a nonzero scalar contained in $Q$, contradicting the fact that $Q$ is a proper ideal.  Therefore $t \not\in Q$, and $\ann_A(V) = A( c - \vartheta)$.  

\begin{rem} {We have seen that  when $\kk$ has  characteristic $0$ and $\phi : R \to R$ is the automorphism of $R = \kk[h,c]$ given by  $\phi(h) = h-1$, $\phi(c) = c$,  
then the $\phi$-stable ideals of $R$ are centrally generated.  The hypotheses of   
Theorem \ref{thm:Whittsimple} are satisfied, and so
$V_{\mathfrak u, \mathcal K} = V_{\mathfrak u}/\mathcal K V_{\mathfrak u}$  is a simple Whittaker module for every maximal ideal   $\mathcal K$  of the center $\mathcal Z$ of $A = R(\phi,t)$.    Thus, when $\kk$ is algebraically closed,    $R^\phi \cap \mathcal K = Q \cap \mathcal Z = R^\phi( c-\vartheta)$ for some
$\vartheta \in \kk$, where $Q = \ann_R(w_{\mathfrak u,\mathcal K}$), and 
$V_{\mathfrak u,\mathcal K} \cong  V_{\mathfrak u}/Q V_{\mathfrak u} = R/Q$, where
$Q = R( c-\vartheta)$ and the action of $A$ is given by
 \eqref{eq:actzero}.
 }
\end{rem}

Now suppose that $\kk$ has characteristic $p > 2$.  Observe that $\phi^p = \hbox{\rm id}_R$ in this case.  Therefore, Proposition \ref{prop:Acent}  implies that the center $\mathcal Z$  of the generalized Weyl algebra $A = R(\phi,t)$ is generated by $X^p, Y^p$ and all the elements of $R$ fixed by $\phi$.   It is clear that $c$ and $h^p - h$ are fixed by $\phi$.   Let us define

$$ z_n = \begin{cases}   h^n  & \quad \hbox{\rm if} \ n \not \equiv 0  \mod p \\
(h^p - h)^{n/p}  & \quad \hbox{\rm if} \ n  \equiv 0  \mod p. \end{cases}$$

\noindent Then $R$ has a basis consisting of the monomials $c^j z_k$ for $j,k \in \z_{\geq 0}$.   Copying the argument of the characteristic $p$ Weyl algebra case (with $t$ replaced by $h$ and the coefficients $g_k$ in that argument assumed to lie in $\kk[c]$ here), we see that  $\mathcal Z$  is generated by $X^p, Y^p, c, h^p-h$.  Observe that $z_p = h^p-h$ and $z_p^j = z_{jp}$.   Thus,  $R^\phi = R \cap \mathcal Z$ is the polynomial algebra  $\kk[c,z_p]$, and $R$ is a free $R^\phi$-module with basis $1,h, \dots, h^{p-1}$. Exactly the same proof as in the Weyl algebra case proves  that $J = R(J \cap \mathcal Z)$ for any $\phi$-stable ideal of $R$. 

Now if $\kk$ is an algebraically closed field of characteristic $p>2$,  and $V$ is a simple Whittaker module, then $Q = \ann_R(w)$ is a maximal $\phi$-stable ideal.  Since $Q \cap \mathcal Z$ is a maximal  ideal of $R^\phi = R \cap \mathcal Z = \kk[c,z_p]$,   there exist $\vartheta, \lambda \in \kk$ so that $Q \cap \mathcal Z = R^\phi\big(c- \vartheta, z_p - (\lambda^p-\lambda)\big)$.  Thus, $Q = R ( c-\vartheta, h^p-h -(\lambda^p-\lambda)\big)$.

The vectors $h^kw$  for $k=0,1,\dots,p-1$ form a basis for the simple module  $V \cong R/Q$.    Since the center of $A$ must act as scalars on $V$, there exists $\alpha \in \kk$ with $X^p = \alpha^p \hbox{\rm id}_V$.   But then $(X-\alpha  \hbox{\rm id}_V)^p = 0$, and the only eigenvalue of $X$ on $V$ is $\alpha$, which must equal $\zeta$.    Since $(h-\lambda)^p - (h-\lambda) =  \prod_{k = 0}^{p-1}\big(h-(\lambda+k)\big)$,  we see that the
vector $v_0: = \prod_{k=1}^{p-1}\big(h-(\lambda+k)\big)w$ satisfies $(h-\lambda)v_0 = 0$.  
Set $v_k = \zeta^{-k} X^k v_0$ and note that $v_k \neq 0$ since $\zeta \neq 0$.    Then the defining relations for $A$ imply that
\begin{eqnarray}\label{eq:paction}  c v_k &=& \vartheta v_k, \qquad 
h v_k = (\lambda+k)v_k     \\
X v_k &=&  \zeta v_{k+1}  \quad (\hbox{\rm subscripts mod} \  p )   \nonumber \\
Y v_k &=&  \frac{1}{2} \zeta^{-1} \Big(\vartheta - r(\lambda +k)\Big)v_{k-1} 
 \quad  (\hbox{\rm subscripts mod} \  p ).  \nonumber
\end{eqnarray}

\noindent  Observe that  $Y^p = \frac{1}{2}  \zeta^{-p}\prod_{k=0}^{p-1} \Big(\vartheta - r(\lambda +k)\Big) \hbox{\rm id}_V$  must hold.    Therefore we have the following:

\begin{thm}\label{thm:smithwhit}  Assume   $\kk$ is algebraically closed of characteristic $p > 2$, and let $A = R(\phi, t)$ be a generalized Weyl algebra over $\kk$ coming from a Smith algebra with defining relations \eqref{eq:smith}, where $s(x)$ is as in \eqref{eq:sfromr}.  If $V$ is a simple Whittaker module for $A$ of type $\zeta$, then $V$ has dimension $p$, and there is a basis $v_k, k = 0,1,\dots, p-1$,  so that the action of $A$ on $V$ is given by \eqref{eq:paction} for scalars $\lambda$, $\vartheta$.  The vector $w = v_0+v_1 + \cdots  + v_{p-1}$ is a cyclic Whittaker vector of type $\zeta$ for $V$.   Then $Q = \ann_R(w) =$  $R\big( c-\vartheta, h^p-h - (\lambda^p-\lambda)\big)$, and $\ann_A(w) = AQ + A(X-\zeta)$.  \end{thm}

\begin{rem} {We have shown that for the Smith algebras of  characteristic $p > 2$  that  the $\phi$-stable ideals of $R = \kk[c,h]$ are centrally generated.  The hypotheses of   
Theorem \ref{thm:Whittsimple} hold, and so
$V_{\mathfrak u, \mathcal K} = V_{\mathfrak u}/\mathcal K V_{\mathfrak u}$  is a simple Whittaker module for every maximal ideal   $\mathcal K$  of $\mathcal Z$.    Thus, when $\kk$ is algebraically closed,    
$V_{\mathfrak u,\mathcal K}$ has dimension $p$ by  Theorem  \ref{thm:smithwhit},
and there is a basis $v_k, k = 0,1,\dots, p-1$,  so that the action of $A$ on $V$ is given by \eqref{eq:paction} for some scalars $\vartheta, \lambda \in \kk$.   The ideal $Q =  \ann_R(w)$ is
generated by the elements   $c-\vartheta, h^p-h - (\lambda^p-\lambda) 
\in Q \cap \mathcal Z = R^\phi \cap \mathcal K$, and 
$V_{\mathfrak u,\mathcal K} \cong  V_{\mathfrak u}/Q V_{\mathfrak u} = R/Q$, where
$Q = R\big( c-\vartheta, h^p-h - (\lambda^p-\lambda)\big)$.  
 }
\end{rem}

\begin{rem}  Suppose $A = R(\phi,t)$ is the Smith algebra defined using the polynomials $s(x) = -1$ and $r(x) = -2x$.  The quotient  $\mathfrak {A}: = A/Ac$ is isomorphic to the Weyl algebra $A_1$, for in $\mathfrak A$ the relation  $YX - XY = 1$ holds and $YX = h+1$.   Let $t' = h+1 = YX$ in $\mathfrak A$.      If   in \eqref{eq:paction}  we set  $\vartheta = 0$ then there is an induced action of $\mathfrak A$ on $V$.   Letting  $\lambda' = \lambda-1$,  we have $t' v_k = (\lambda' + k) v_k,  \quad X v_k = \zeta v_{k+1}, \quad Yv_k = \zeta^{-1}(\lambda'+k-1)v_k$, which are precisely the relations we obtained in Theorem \ref{thm:Weylp} for the Whittaker modules of a Weyl algebra in characteristic $p$. \end{rem}

\section{Quantum Smith algebras} \label{sec:qSmith}

In \cite{JiWangZhou:quantumSmith} Ji, Wang, and Zhou introduced a family of associative algebras $\mathcal A$ which generalize the quantized  enveloping algebra $U_q(\mathfrak{sl}_2)$ and which are quantum versions of the Smith algebras  with defining polynomial $s(h) = h^{m+1} - h^m$ for some $m$.     When the underlying field is the complex numbers, Tang \cite{tang:WMSmithAlg} determined the irreducible weight modules for these algebras, showed that the finite-dimensional modules for $\mathcal A$ are weight modules  which are completely reducible, and  obtained analogues of the results by Ondrus in \cite{ondrus:Whittaker} for  the Whittaker $\mathcal A$-modules.

The algebras $\mathcal A$ have a realization as generalized Weyl algebras, and here we illustrate how  results we have obtained can be specialized to recover  results in \cite{tang:WMSmithAlg} and \cite{ondrus:Whittaker}.  We also determine all the simple Whittaker modules in the root of unity case, which is not considered either in \cite{tang:WMSmithAlg} or in \cite{ondrus:Whittaker}.

Let $\kk$ be a field of characteristic not 2 and fix an integer $m \geq 1$.  Assume  $q \in \kk, q \neq 0,\pm 1,$ and $q^2$ is not an $m$th root of unity.   Consider a unital associative algebra $\mathcal A$ over $\kk$ with generators $E,F,K^{\pm 1}$ which satisfy  the defining relations
\begin{eqnarray}\label{eq:qsmith} &&KE = q^2 EK,  \qquad  KF = q^{-2}FK \qquad  KK^{-1} = 1 = K^{-1}K \\   && EF-FE = \frac{K^{m}-K^{-m}}{q-q^{-1}}, \nonumber \end{eqnarray} where $m \in \z_{\geq 1}$.   In particular, when  $m = 1$, the algebra $\mathcal A$ is isomorphic to  $U_q(\mathfrak{sl}_2)$.  Tang gave a realization of this algebra as a hyperbolic algebra (as defined in \cite{rosenberg:quantizedAlg}).  Here we realize it as a generalized Weyl algebra (the two realizations are equivalent).  The element 
$$c = FE + \frac{q^mK^{m} + q^{-m}K^{-m}}{(q^m-q^{-m})(q-q^{-1})}$$ 
is central,  and it generates the center of  $\mathcal A$ when $q^2$ is not a root of unity.  (See  \cite[Lem. 3.1.1 and Prop. 3.1.2]{tang:WMSmithAlg}.)

Let $R = \kk[K^{\pm 1}, c]$  and define an automorphism $\phi$ on $R$ by setting $\phi(K^{\pm 1}) = q^{\mp 2}K^{\pm 1}$ and $\phi(c) = c$.   Let

\begin{equation} \label{eq:smithqt}  t = c- \frac{q^mK^{m} + q^{-m}K^{-m}}{(q^m-q^{-m})(q-q^{-1})}, \end{equation}
and assume $A = R(\phi, t)$ is the generalized Weyl algebra constructed from this data.   Thus in $A$  we have $YX = t$,  and
$$XY = \phi(t) = c- \frac{q^{-m}K^{m} + q^{m}K^{-m}}{(q^m-q^{-m})(q-q^{-1})},$$
and $A$ can be seen to be isomorphic to $\mathcal A$ by identifying $X$ with $E$ and $Y$ with $F$.

We begin by describing the center $\mathcal Z$ of $A$.  If $q^2$ is not a root of unity, then $\mathcal Z = \kk[c]$ by Proposition \ref{prop:Acent}.   If $q^2$ is a primitive $\ell$th root of unity for $\ell \neq m$, then $\phi^\ell = \hbox{\rm id}_R$.    Therefore, it follows from Proposition \ref{prop:Acent} that the center of the generalized Weyl algebra $A = R(\phi,t)$ is generated by $X^\ell, Y^\ell$ and all the elements of $R^\phi$.  It is clear that $\phi(K^{\pm \ell}) = K^{\pm \ell}$.    Suppose $h = \sum_{j=-r}^s  h_j(c) K^j$ is fixed by $\phi$.  Then  for each $j$ with $h_j(c) \neq 0$, we must have $q^{-2j} = 1$, or $j \equiv 0 \mod \ell$.   Thus in this case,  the center $\mathcal Z$  of $A$ is generated by $X^\ell, Y^\ell, c, K^{\pm \ell}$.

Let $\Xi = \{q^{2k} \mid k \in \z \}$, and note that the algebra $R$ decomposes into eigenspaces $R_\xi$, $\xi \in \Xi$,  relative to $\phi$, where $R_\xi = \{ r\in R \mid \phi(r) = \xi r\}$.   Thus, 
$$R_\xi = {\rm span}_\kk \setof{c^iK^j}{\mbox{$i \in \z_{\ge 0}$ and $q^{-2j} = \xi$}}$$ 
and $R = \bigoplus_{\xi} R_\xi$ gives a grading of $R$.   If $J$ is a $\phi$-stable ideal  of $R$, then $J = \bigoplus_{\xi} J_\xi$ where $J_\xi = J \cap R_\xi$.   For $\xi = q^{-2n} \in \Xi$, 
$$K^{-n}J_\xi \subseteq \left\{ \begin{array}{ll} J \cap \kk[c] & \mbox{if $q^2$ is not a root of unity} \\ J \cap \kk[c,K^{\pm \ell}] & \mbox{if $q^2$ is a primitive $\ell$th root of unity,} \end{array} \right.$$
so by the previous paragraph, we have shown that $K^{-n}J_\xi \subseteq J \cap \mathcal Z$ whether or not $q^2$ is a root of unity.  This implies the following.

\begin{lem}\label{lem:quantcentgen}
Let $J$ be a $\phi$-stable ideal of $R = \kk[c,K^{\pm 1}]$, where $\phi (c) = c$ and $\phi (K) = q^{-2}K$.  Then $J = R(J \cap \mathcal Z)$, where $\mathcal Z$ is the center of the corresponding generalized Weyl algebra $A = R(\phi,t)$ with $t$ as in \eqref{eq:smithqt}.
\end{lem}

Now let $(V,w)$ be a Whittaker pair for $A$ of type $\zeta$ with $Q = \ann_R(w)$.   As $Q$ is $\phi$-stable, $Q$ is centrally generated and
\begin{equation}\label{annrwsmith} 
Q = \ann_R(V) = \left\{ \begin{array}{ll} R \mathcal Z_V & \mbox{if $q^2$ is not a root of unity} \\ R(\mathcal Z_V \cap R) & \mbox{if $q^2$ is a root of unity,} \end{array} \right.  \end{equation} 
\noindent where $\mathcal Z_V = \ann_A(V) \cap \mathcal Z$.  Note that $A \mathcal Z_V \subseteq \ann_A(w)$ regardless of whether $q^2$ is a root of unity.  Thus by Theorem \ref{thm:ann_A(w)}, 
\begin{equation}\label{eq:annwsmith} \ann_A(w) = AQ + A(X-\zeta)  =  A \mathcal Z_V  + A (X-\zeta).
\end{equation}
\noindent   When $q^2$ is not a root of unity, this is Theorem 3.2.2 of \cite{tang:WMSmithAlg}.  Theorem 3.2.3 of \cite{tang:WMSmithAlg} establishes a one-to-one correspondence between isomorphism classes of Whitaker modules for $A$ of type $\zeta$ and ideals of the center $\mathcal Z = \kk[c]$ given by $V \rightarrow \mathcal Z_V$.   In the present setting (with no assumption on $q^2$), Theorem \ref{thm:bijectIdealsModules} gives a bijection $V \mapsto \ann_R(w) = Q$ between  isomorphism classes of Whitaker modules for $A$ of type $\zeta$  and $\phi$-stable ideals of  $R$.  But there is a bijection between $\phi$-stable ideals $J$ and ideals of $\mathcal Z$ given by $J \mapsto J \cap \mathcal Z$ since $J = R(J \cap \mathcal Z)$.

Assume that $q^2$ is not a root of unity, and let $V$ be a simple Whittaker module, with $\kk$ algebraicially closed.  Since $Q := \ann_R(w)$ is generated by its intersection with $\kk[c]$, there must exist $\vartheta \in \kk$ such that $Q = R(c - \vartheta)$.  Notice that $R/Q \cong \kk[K^{\pm 1}]$ is a domain, and it is clear that the induced automorphism $\overline \phi: R/Q \to R/Q$ has infinite order since $\overline \phi (K) = q^{-2}K$.  Thus as long as $t \not\in Q$, Theorem \ref{thm:annVwhenQprime} implies that $\ann_A(V) = A(c - \vartheta)$.  Recall that $t = c- \frac{q^mK^{m} + q^{-m}K^{-m}}{(q^m-q^{-m})(q-q^{-1})}$.  If $t \in Q$, then $(\phi - {\rm id}_R)(\phi - q^{2m}{\rm id}_R)(t)$ is a nonzero multiple of $K^m$ belonging to $Q$, contradicting the fact that $Q$ is a proper ideal.  Therefore $t \not\in Q$, and $\ann_A(V) = A(c - \vartheta)$.

Suppose now that $q^2$ is a primitive $\ell$th root of unity for $\ell \neq m$.    Assume $\kk$ is algebraically closed  and let $V$ be a simple Whittaker module for $A$ with Whittaker vector $w$ of type $\zeta$ and with  $Q = \ann_R(w)$.    Then  $Q \cap \mathcal Z$ is a maximal ideal of $\kk[c,K^{\pm \ell}]$, and there exist scalars $\vartheta, \lambda$,  with $\lambda \neq 0$, so that $Q \cap \mathcal Z = R^\phi(c-\vartheta, K^{\pm \ell}-\lambda^{\pm \ell})$,  and $Q = R(c-\vartheta, K^{\pm \ell}-\lambda^{\pm \ell})$.    Thus, $V \cong R/Q$ has a basis consisting of the vectors  $K^j w$ for $j=0,1, \dots, \ell-1$.    Since $X K^j w = q^{-2j}\zeta K^j w$,   we see that $X^\ell = \zeta^\ell \hbox{\rm id}_V$.  If $v_0 := \sum_{j=0}^{\ell-1}  \lambda^{-j} K^j w$, then $K v_0 = \lambda v_0$.     The vectors $v_j = \zeta^{-j} X^j v_0$ for $j = 0,1,\dots, \ell-1$ are eigenvectors for $K$ ($K v_j = \lambda q^{2j}v_j$)  corresponding to different eigenvalues.  Hence they are linearly independent and comprise a basis for $V$.   The action of $A$ relative to this basis is given by 
 
\begin{eqnarray}\label{eq:laction}  c v_j &=& \vartheta v_j. \qquad \quad 
K v_j =  \lambda q^{2j} v_j    \\
X v_j &=&  \zeta v_{j+1}  \quad (\hbox{\rm subscripts mod} \  \ell )   \nonumber \\
Y v_j &=&  \zeta^{-1} 
\Big(\vartheta - \frac{\lambda^m q^{(2j+1)m} + \lambda^{-m}q^{-(2j+1)m}}{(q^m-q^{-m})(q-q^{-1})}\Big)v_{j-1} 
 \quad  (\hbox{\rm subscripts mod} \  \ell ).  \nonumber
\end{eqnarray}  
\noindent     Note that   
$Y^\ell =  \displaystyle{\zeta^{-\ell} \prod_{j=0}^{\ell-1}
\Big(\vartheta - \frac{\lambda^m q^{(2j+1)m} + \lambda^{-m}q^{-(2j+1)m}}{(q^m-q^{-m})(q-q^{-1})}\Big)
\hbox{\rm id}_V}$  must hold.    Therefore, we have the following.

\begin{thm}\label{thm:qsroot1} Assume   $\kk$ is algebraically closed of characteristic $\neq 2$, and let $A = R(\phi, t)$ be a generalized Weyl algebra over $\kk$ coming from a quantum Smith algebra with defining relations \eqref{eq:qsmith},  where $q^2$ is a primitive $\ell$th root of unity, and $\ell \neq m$.   If $V$ is a simple Whittaker module for $A$ of type $\zeta$, then $V$ has dimension $\ell$, and there is a basis $v_j, j = 0,1,\dots, \ell-1$,  so that the action of $A$ on $V$ is given by \eqref{eq:laction} for scalars $\lambda$, $\vartheta$.   The vector $w = v_0+v_1 + \cdots  + v_{\ell-1}$ is a cyclic Whittaker vector of type $\zeta$ for $V$.   Then $Q = \ann_R(w) = R\big( c-\vartheta, K^{\pm \ell} - \lambda^{\pm \ell}\big)$,  and $\ann_A(w) = AQ + A(X-\zeta)$.    \end{thm}

\begin{rem} {Since $R^\phi$ is $\kk[c]$ (if $q^2$ is not a root of unity) or $\kk[c,K^{\pm \ell}]$ (if $q^2$ is a primitive $\ell$th root of unity), it follows that $R^\phi$ is always finitely generated over $\kk$.  Lemma \ref{lem:quantcentgen} says that every $\phi$-stable ideal of $R$ is centrally generated.  Thus Theorem \ref{thm:Whittsimple} implies that $V_{\mathfrak u, \mathcal K}$ is a simple Whittaker module for every maximal ideal $\mathcal K$ of 
$\mathcal Z$.  When $q^2$ is not a root of unity, then there is a scalar $\vartheta$ so that
$V_{\mathfrak u,\mathcal K} \cong V_{\mathfrak u}/Q V_{\mathfrak u} = R/Q$, where
$Q = R(c-\vartheta)$.    When $q^2$ is a primitive $\ell$th root of unity,  
$V_{\mathfrak u, \mathcal K}$ has dimension $\ell$ by  Theorem \ref{thm:qsroot1},
and the
action of $A$ on $V_{\mathfrak u, \mathcal K}$ is given by \eqref{eq:laction} for scalars $\lambda$, $\vartheta$ such that $c-\vartheta$, $K^{\pm \ell} -\lambda^{\pm \ell}\in Q \cap \mathcal Z = R^\phi \cap \mathcal K$.   Thus, $V_{\mathfrak u,\mathcal K} \cong V_{\mathfrak u}/Q V_{\mathfrak u}
= R/Q$, where $Q = R\big( c-\vartheta, K^{\pm \ell} -\lambda^{\pm \ell}\big)$.  }
\end{rem}

\noindent  Department of Mathematics, University
of Wisconsin, Madison, Wisconsin 53706, USA,  
\  \textbf{e-mail:} \ {\tt benkart@math.wisc.edu}  
\medskip
 
\noindent Department of Mathematics,
 Weber State University, Ogden, UT  84408, USA,   
\  \textbf{e-mail:} \ {\tt MattOndrus@weber.edu}

\end{document}